\documentclass[review,onefignum,onetabnum]{siamart250211}
\usepackage{tikz,subfigure}
\usetikzlibrary{shapes.geometric}
\usepackage{booktabs}
\usepackage{tabularx}
\usepackage{array}
\newcolumntype{C}[1]{>{\centering\let\newline\\\arraybackslash\hspace{0pt}}m{#1}}
\usepackage{algorithmic}
\usepackage{graphicx}
\usepackage{amssymb}
\newcommand{\grad}{\nabla}
\newcommand{\REMOVE}[1]{}

\newcommand{\half}{\frac{1}{2}}
\newcommand{\StencilTwoD}[4]{ %
  \left( %
    \begin{array}{ccccc}%
      #4 &  & #2 & & #4\\%
      &&&&\\%
      #3 &  & #1 & & #3\\%
      &&&&\\%
      #4 &  & #2 & & #4\\%
    \end{array}\right)}

\newcommand{\StencilThreeD}[9]{ %
\left\{
\begin{array}{cl}
  #9 \StencilTwoD{#1}{#3}{#2}{#5}, & \stackrel{\mbox{middle}}{\mbox{xy-plane}}\\
&\\[10pt]
  #9 \StencilTwoD{#4}{#7}{#6}{#8}, & \stackrel{\mbox{top/bottom}}{\mbox{xy-plane}}
  \end{array}
\right.
}

\newcommand{\CMS}[1]{\textcolor{red}{\textbf{CMS: #1}}}

\newcommand{\RST}[1]{\textcolor{green}{\textbf{RST: #1}}}
\newcommand{\squots}{\textquotesingle \textquotesingle~~~~~~}
\definecolor{color1}{RGB}{255,174,120}
\definecolor{purpleheart}{rgb}{0.41, 0.21, 0.61}
\definecolor{yellow(ncs)}{rgb}{1.0, 0.83, 0.0}
\definecolor{yaleblue}{rgb}{0.06, 0.3, 0.57}
\definecolor{venetianred}{rgb}{0.78, 0.03, 0.08}
\definecolor{tealgreen}{rgb}{0.0, 0.51, 0.5}

\usepackage{lipsum}
\usepackage{amsfonts}
\usepackage{graphicx}
\usepackage{epstopdf}
\usepackage{algorithmic}
\ifpdf
  \DeclareGraphicsExtensions{.eps,.pdf,.png,.jpg}
\else
  \DeclareGraphicsExtensions{.eps}
\fi

\newsiamremark{remark}{Remark}
\newsiamremark{hypothesis}{Hypothesis}
\crefname{hypothesis}{Hypothesis}{Hypotheses}
\newsiamthm{claim}{Claim}

\headers{Improving Smoothed Aggregation Robustness}{Siefert, Tuminaro, Sunderland}

\title{Improving Smoothed Aggregation AMG Robustness on Stretched Mesh Applications\thanks{Submitted to the editors 2/2026.
\funding{This work was partially supported by the Laboratory Directed Research and Development program (Project 236939). Tuminaro was also partially supported by the U.S.~Department of Energy, Office of Science, Office of Advanced Scientific Computing Research, Applied Mathematics program.  Sandia National Laboratories is a multimission laboratory managed and operated by National Technology and Engineering Solutions of Sandia, LLC., a wholly owned subsidiary of Honeywell International, Inc., for the U.S. Department of Energy's National Nuclear Security Administration under grant~DE-NA-0003525.  This paper describes objective technical results and analysis.  Any subjective views or opinions that might be expressed in the paper do not necessarily represent the views of the U.S. Department of Energy or the United States Government.}}} 

\author{
  Christopher M. Siefert\thanks{Scalable Algorithms Department, Sandia National Laboratories, Albuquerque, NM}
  (\email{csiefer@sandia.gov}, \url{https://www.sandia.gov/-csiefer/staff/chris-siefert})
\and
  Raymond S. Tuminaro\thanks{Computational Mathematics Department, Sandia National Laboratories, Livermore, CA}
  (\email{rstumin@sandia.gov}, \url{https://www.sandia.gov/ccr/staff/raymond-s-tuminaro})
\and
  Daniel J. Sunderland\thanks{Sunderland Development, Salem, Oregon.\email{dan@sunderland.site}}
  }

\usepackage{amsopn}

\ifpdf
\hypersetup{
  pdftitle={Distance Laplacian Strength-of-Connection and Thresholding Methods for Smoothed Aggregation Multigrid}
  pdfauthor={C. Siefert, R. Tuminaro, and D. Sunderland}
}
\fi

\begin{document}
\nolinenumbers

\maketitle

\begin{abstract}
Strength-of-connection algorithms play an essential role 
within algebraic multigrid (AMG).
Specifically, the strength-of-connection scheme determines which matrix nonzeros are classified
as {\it weak} and so ignored
when coarsening matrix graphs and defining sparsity patterns for interpolation.
The general goal is to 
encourage coarsening only in directions
where error can be smoothed 
and to also avoid coarsening
across sharp problem variations. 
Unfortunately, developing robust, practical, and inexpensive strength-of-connection schemes is
extremely challenging.

The classification of matrix nonzeros involves four distinct aspects:
(a) choosing a strength-of-connection matrix,
(b) scaling its values, 
(c) choosing a criterion to classify the scaled values
as either {\it strong} or {\it weak}, and
(d) dropping weak entries which includes adjusting matrix values to account for dropped terms.
Typically, smoothed aggregation AMG uses the linear system being solved as the strength-of-connection matrix.
It scales these values symmetrically using the square-root of the matrix diagonal. It classifies
based on whether scaled values are above or below a user-supplied threshold. Finally,
it adjusts matrix values by modifying the diagonal so that the sum of entries within each row of the
dropped matrix matches the sum of the original. While these procedures often work well, we illustrate
stretched mesh failure cases that motivate
alternatives to improve robustness. The first alternative uses a
distance Laplacian strength-of-connection matrix. The second idea centers on
non-symmetric scaling algorithms. 
We then investigate alternative classification criteria
based on
identifying a significant gap in the values of the scaled entries.
Finally, an alternative lumping procedure is proposed where row sums are preserved by
modifying all retained matrix entries (as opposed to just diagonal entries).
A series of numerical results is given to illustrate algorithm trade-offs demonstrating
in some cases notably more robust convergence
on matrices coming from  linear finite elements on stretched meshes.
\end{abstract}



\section{Introduction}\label{sec:intro}
Algebraic multigrid (AMG) \cite{BrMcRu84,RuSt85} is a popular
solver for large, sparse linear systems. 
The main idea is to employ coarse problems to accelerate convergence
by leveraging the fact
\REMOVE{
Simple
iterative relaxation methods are then used to smooth error on the various
meshes.
Classical AMG
\cite{BrMcRu84,RuSt85} generates coarse meshes through the careful
selection of a subset of the fine grid unknowns to also be represented
on the coarse grid.  Interpolation weights are then calculated based
on neighborhood relations between fine and coarse unknowns.}
%
that simple relaxation schemes effectively damp error components whose
$A$-norm is {\it relatively} large\footnote{The $A$-norm of a vector $v$ is given by $\sqrt{v^T A v}$ where $A$ is the positive-definite symmetric matrix that one is interested in solving.}, and so errors that remain after
relaxation have relatively small $A$-norm.  Ideally, these small $A$-norm
components are smooth and can be well approximated on a lower resolution mesh. Thus,
one projects a correction equation to the next coarsest mesh and recursively repeats the process
of relaxation (damping oscillatory errors
relative to the mesh resolution) followed by further projection.
In many cases, however, some oscillatory error components
may have small $A$-norm that is difficult to damp.
As is well known (e.g., see \cite{savariants22}), this issue arises with either anisotropic
partial differential equations (PDEs)
or isotropic PDEs that are discretized on stretched meshes. Intuitively,
most relaxation schemes primarily smooth error in dominant directions.
In the mesh-stretching case, these would be directions where grid points are
relatively close to each other.
To 
avoid convergence degradation,
either special relaxation schemes must be devised (e.g., line relaxation on structured grids)
to damp some oscillatory error components
or instead 
coarsening must not occur in directions where errors are oscillatory after relaxation.
Most AMG algorithms adopt the later strategy, which we also target in this paper.

While anisotropic challenges have been studied,
a robust multigrid algorithm
remains elusive
for applications that employ unstructured meshes.
Earlier anisotropic multigrid research focused on structured grid approaches
(e.g.~\cite{semijones,De1982, DeIdRu1992,Sc1989}) usually employing semi-coarsening ideas
sometimes in conjunction with line/plane relaxation. However, for more
general meshes there are usually no lines or planes to leverage
(see~\cite{MAVRIPLIS1998141} for some unstructured adaptations).
Overall,  the most significant unstructured mesh impediment is that of detecting the
directions where coarsening should be restricted, as doing this robustly
can be computationally expensive.  Most
techniques~\cite{ea0982c56c9743fcadcc133e274ae803,8139e1fa7c69487785b0685345067e5c,BrBrMaMaMc06,doi:10.1137/17M1123456,BrZi2007,BrannickF10,
Livne04,OlScTu10}
employ some kind of analysis phase to discover directions where relaxation struggles to smooth errors.
Not surprisingly, the more robust techniques are also generally
fairly costly. Instead, most popular AMG codes rely on simple
heuristics that attempt to inexpensively infer coarsening directions from matrix coefficients.
While these heuristics frequently work, they fail all too often. In this paper,
we illustrate cases where smoothed aggregation AMG (SA)~\cite{VaBrMa01,VaMaBr96} fails in the
context of stretched meshes.  Our focus on stretched meshes (as opposed to anisotropic PDEs)
is motivated by applications where anisotropic behavior frequently arises from
stretched meshes. 
Further, we concentrate on matrices produced via first-order finite elements as these
are heavily used by scientists. Finally, while we examine
Poisson problems, many of our example meshes come from complex applications
(e.g., fluid flow or Maxwell's equations), which commonly rely on Poisson-like sub-solves or
which give rise to similar issues.

In the smoothed aggregation context, the standard heuristics for anisotropic applications
rely on strength-of-connection
ideas that were originally introduced for classical AMG~\cite{BrMcRu84,RuSt85}.
The basic strategy encourages coarsening in preferential directions by dropping
entries from the matrix used within the algorithm's coarsening
and interpolation construction phases.  In this paper, we take a somewhat different
perspective on the dropping procedure in that we highlight the four independent
individual sub-steps that together determines the process.
The first defines a strength-of-connection (SOC) matrix that is used in the classification process.
Normally, this is taken as the original system $A$ that is being solved, but one could instead
define a matrix with the same sparsity pattern as $A$ but different values for the entries.
The second sub-step scales the SOC matrix entries so that each scaled entry
reflects a relative size as opposed to the absolute sizes of the unscaled system. In smoothed
aggregation a symmetric scaling employing the matrix diagonal is normally used.
The third sub-step applies a criterion that examines the scaled entries and determines
whether this entry should be kept or dropped (because it is strong or weak).
The fourth sub-step adjusts the dropped matrix so that it retains some properties
of the original matrix.  The standard
smoothed aggregation procedure simply alters the matrix diagonal of the dropped matrix so that
the row sum (i.e., sum of entries within each row) of the dropped matrix matches the corresponding
row sum of the original.

The heart of the paper discusses the above four smoothed aggregation sub-steps 
highlighting deficiencies with the standard approach by first examining a simple Poisson operator
on a structured mesh that is stretched in one direction. For the first sub-step, we
evaluate the use of a distance Laplacian system as an alternative to define the
SOC matrix. For the second sub-step we consider algorithms that scale the SOC matrix entries in a non-symmetric
fashion. This includes classical AMG~\cite{BrMcRu84,RuSt85} as well as some alternatives that we propose
such as the {\it cut-drop} algorithm. For the third sub-step we consider an approach which
seeks to find a gap between scaled entries within each matrix row, classifying those
on one side of the gap as strong and those on the other side as weak.
For the final sub-step,
we propose an alternative distributed lumping procedure {\it distrib-lump} which
preserves the row sum without changing the sign of the entries along the
matrix diagonal.
We compare and contrast the different sub-step algorithms demonstrating the utility of the
best-performing algorithms on a set of
unstructured test cases,
illustrating how the new approaches can
avoid some standard SA pitfalls  that arise from stretched meshes. 

\section{Smoothed Aggregation and its Shortcomings}\label{sec:dropping_survey}
In this section, smoothed aggregation AMG is described focusing on the four
sub-steps that are intended to address anisotropic phenomena.
Before doing this, however, we present a structured stretched model problem that will be used to highlight
smoothed aggregation shortcomings by contrasting it with
a proper  structured geometric multigrid algorithm.

\subsection{Structured Geometric Multigrid Semi-Coarsening}

As noted, mesh coarsening should only occur in directions where relaxation is effective.
To illustrate this, consider a Poisson problem
\begin{eqnarray}\label{eq:poisson}
-\Delta u = g &&              
\end{eqnarray}
on a cuboid domain $\Omega$ with Dirichlet conditions on all boundaries except the two $yz$-planes, which
have Neumann boundary conditions.  The discrete linear system
\begin{equation} \label{eq: model linear system}
A u = f
\end{equation}
is obtained by employing tri-linear hexahedral nodal finite elements on a
tensor-product mesh with constant mesh spacing in the $x$
and $y$ directions given by $h$. The constant mesh spacing in the $z$ direction is given
by $\alpha h$.

Consider now a semi-coarsening structured geometric multigrid method, $\mbox{\sc semi}_{\bar{\alpha}}$, that always coarsens
in the $x$ and $y$ directions but only coarsens in the $z$ direction if the $z$ mesh spacing
is less than or equal to $\bar{\alpha}$ times the $x$ spacing on that mesh. Thus, $\mbox{\sc semi}_{\infty}$ corresponds
to standard multigrid that coarsens in all directions. When $\alpha \gg 1$, a  $\mbox{\sc semi}_{1}$ algorithm
would initially construct coarse meshes by only coarsening in the $x$ and $y$ directions. Should this
lead to a coarse mesh where the mesh spacing is identical in all directions, then it would coarsen in all directions
to generate additional coarse meshes.  When coarsening
occurs, we always consider factors of three along each axis (which is a typical SA coarsening rate in one dimension).

Table~\ref{tab:regular 3D meshes} shows iteration counts
\begin{table}[h!]
\centering
\caption{CG iterations ($10^{-10}$ residual reduction tolerance) with $4$ level multigrid  V-cycle preconditioner and 1 pre- and 1 post-Jacobi relaxation sweep on all levels, except coarsest which uses a direct solver. The Jacobi damping parameter is $.6$. Quote marks indicate identical coarsening as previous line.}
\label{tab:regular 3D meshes}
\begin{tabular}{c|lccccc}
                           &         &  ~~~~~~       &  ~~~~~~ &   ~~~~~~ &   ~~~~~~ &   ~~~~~~ \\[-10pt]
\toprule
coarsen                                       \\
algorithm \\[-16pt]
                           &$\alpha=$&        1.     &      3.     &      9.      & 27.      & 81.   \\[5pt]
\midrule
$\mbox{\sc semi}_1 $       &         &     17        &     17      &     22       & 23       & 23    \\
$\mbox{\sc semi}_3 $       &         &~~~~~\squots   &     43      &     40       & 30       &~~~~~\squots\\
$\mbox{\sc semi}_9 $       &         &~~~~~\squots   &~~~~~\squots &     120      & 86       & 44    \\
$\mbox{\sc semi}_{27}$       &         &~~~~~\squots   &~~~~~\squots &~~~~~\squots  & 266      & 127   \\
$\mbox{\sc semi}_\infty$   &         &~~~~~\squots   &~~~~~\squots &~~~~~\squots  & ~~~~~\squots  & 394   \\
    \bottomrule
\end{tabular}
\end{table}
Jacobi-smoothed  $\mbox{\sc semi}_{\bar{\alpha}}$ is applied to a $82 \times 82 \times  82 $ problem.
The finest mesh spacing in the $x$ and $y$ direction is fixed at $1/81$ and so the computational domain might be more elongated
in the $z$ direction depending on the value of $\alpha$. As an example, a $\mbox{\sc semi}_{3}$ algorithm applied to
a fine mesh with $\alpha = 9$ would use the following mesh sequence:
$ 82  \times 82  \times 82 $, $ 28 \times 28 \times 82$,   $ 10 \times 10 \times 28 $, and $ 4 \times 4 \times 10$.
That is, coarsening in the $z$ direction only starts on the $ 28 \times 28 \times 82$ mesh as at this point the $z$ direction spacing
is three times larger than that of the $x$ and $y$ directions.
Multigrid interpolation employs tri-linear interpolation,
restriction is the transpose of interpolation.
It is clear that $\mbox{\sc semi}_{1}$ is generally
the best performing algorithm in terms of iteration counts as it does not coarsen in the $z$ direction until its mesh spacing is
comparable to that of the $x$ and $y$ directions.
Clearly, full isotropic coarsening ($\mbox{\sc semi}_{\infty}$) struggles
when the $z$ direction mesh spacing is large relative to the $x$ and $y$ spacing.
In this scenario, damped-Jacobi
$$
u^{(k+1)} \leftarrow u^{(k)}  + \omega D^{-1} (f - A u^{(k)} )
$$
only smooths errors in the $x$ and $y$ directions.
Here, $u^{(k)}$ refers to the approximate solution after $k$ Jacobi iterations,
$\omega$ is a suitable damping parameter, and $D$ is the diagonal of $A$.  We note that in these experiments the cost per iteration
(not shown here) is slightly higher when one compares $\mbox{\sc semi}_{1}$ to $\mbox{\sc semi}_{\infty}$, but the iteration
savings far outweigh this additional cost. Using this as a guide, we can see that an algebraic method when applied to
an isotropic diffusion equation should avoid
coarsening in directions associated with neighboring grid points that are relatively far from the central point.
We now describe smooth aggregation and investigate how it behaves on the above model problem.

\subsection{Smoothed Aggregation} \label{sec: sa}
The first step in applying smoothed aggregation to \eqref{eq: model linear system} is to coarsen the graph associated with $A$.
In the simplest case, one can apply an aggregation algorithm directly to the graph of the matrix $A$ with the goal of assigning
vertices to each aggregate, ${\cal A}_k$, so that the resulting ${\cal A}_k$'s are disjoint, and approximately the same
size. These aggregates are effectively coarse mesh points that will be used to define the next mesh in the multigrid
hierarchy.
A common approach is to choose a vertex that has not yet been assigned to any aggregate and is not a neighbor of
an already created aggregate (i.e., no edge connects this vertex to an already assigned vertex). A new initial aggregate is then formed
by grouping this chosen vertex with its neighboring vertices. This process is repeated to continue aggregation until only unassigned
vertices remain that have at least one neighbor that is assigned to a created aggregate. Additional heuristics are then applied to either
merge unassigned vertices into existing aggregates or create a new aggregate by grouping together a cluster of unassigned
vertices~\cite{VaMaBr96,sc2000}.

While aggregation might be applied directly to $A$, we can mimic semi-coarsening if we remove some of $A$'s entries.
Suppose that for our example problem an auxiliary $\tilde{A}$ matrix is defined by removing nonzero $A_{ij}$ when the $z$
coordinate of the $i^{th}$
and $j^{th}$ mesh points satisfy $z_i \ne z_j$. In this case, each aggregate will lie entirely within a single $xy$-plane.
Thus no coarsening is performed in the $z$-direction, corresponding to a coarsening scheme that resembles the
advantageous $\mbox{\sc semi}_{1}$ method when $ \alpha \gg 1$. Of course, if we do not know that the $z$ direction
corresponds to the distant points, then this must somehow be detected. A traditional smoothed aggregation method will attempt
to infer this based on whether scaled matrix entries are above or below a user-supplied threshold $\theta\in[0,1]$.
Specifically, smoothed aggregation~\cite{VaBrMa01,VaMaBr96} employs the following criterion for retaining
$A_{ij}$ entries:
\begin{equation}\label{eq:crit_sa}
|A_{ij}| \geq \theta \sqrt{A_{ii} A_{jj}} .
\end{equation}
Following classical AMG, these retained entries are referred to as ``strong'' while discarded entries are
termed ``weak''.  One advantageous aspect of this criterion is that it is symmetric in nature. That is, if
$A$ is symmetric and $A_{ij}$ satisfies \eqref{eq:crit_sa}, then $A_{ji}$ must also satisfy \eqref{eq:crit_sa}
implying that $\tilde{A}$ retains $A$'s symmetry property. In contrast to smoothed aggregation, classical AMG~\cite{BrMcRu84,RuSt85}
employs a different criterion
\begin{equation} \label{eq: tradclass}
-A_{ij} \geq \theta \max_{k\not=i} -A_{ik} .
\end{equation}
Though this criterion is not symmetric, it does have advantages (discussed shortly) over \eqref{eq:crit_sa}.
Notice that this criterion drops all positive off-diagonals entries if at least one off-diagonal entry is negative.
While this generally works fine for standard discretizations of Poisson operators, it may not be
valid for more sophisticated PDEs or for exotic discretizations as it assumes that any positive entry must be weak.

Figure~\ref{fig:bleez}
\begin{figure}[htb!]\label{fig:bleez}
\begin{equation}\label{eq:sten3d_A}
A_{i:} =
\StencilThreeD{16+32\alpha^2}
{4-4\alpha^2}
{4-4\alpha^2}
{-8+8\alpha^2}
{1-4\alpha^2} 
{-2-\alpha^2}
{-2-\alpha^2}
{-\half-\alpha^2}
{\frac{h}{18\alpha}}
\end{equation}
\vskip -2.16in
\hskip 1.46in {\begin{tikzpicture}
\draw[dashed,tealgreen,line width=2pt] (3.51,0.18) rectangle (5.10,0.83);
\draw[purpleheart,line width=2pt] (6.36,0.51) node [ellipse,draw,minimum width=.7in,outer sep=-10pt, minimum height=.21in]{};
\draw[dashed,tealgreen,line width=2pt] (7.54,0.18) rectangle (9.10,0.83);
\draw[purpleheart,line width=2pt] (4.31,1.36) node [ellipse,draw,minimum width=.7in,outer sep=-10pt, minimum height=.21in]{};
\draw[venetianred,dotted,line width=2pt] (6.35,1.36) node [ellipse,draw,minimum width=.7in,outer sep=-10pt, minimum height=.21in]{};
\draw[purpleheart,line width=2pt] (8.37,1.36) node [ellipse,draw,minimum width=.7in,outer sep=-10pt, minimum height=.21in]{};
\draw[dashed,tealgreen,line width=2pt] (3.51,1.86) rectangle (5.10,2.51);
\draw[purpleheart,line width=2pt] (6.31,2.20) node [ellipse,draw,minimum width=.7in,outer sep=-10pt, minimum height=.21in]{};
\draw[dashed,tealgreen,line width=2pt] (7.54,1.86) rectangle (9.10,2.51);
\draw[yellow(ncs),line width=2pt] (3.55,3.13) rectangle (5.10,3.72);
\draw[dashed,line width=2pt,yaleblue] (5.55,3.13) rectangle (7.10,3.72);
\draw[yellow(ncs),line width=2pt               ] (7.55,3.13) rectangle (9.10,3.72);

\draw[dashed,line width=2pt,yaleblue] (3.55,3.92) rectangle (5.10,4.51);
\draw[dashed,line width=2pt,yaleblue] (7.55,3.92) rectangle (9.10,4.51);
\draw[yellow(ncs),line width=2pt               ] (3.55,4.71) rectangle (5.10,5.32);
\draw[dashed,line width=2pt,yaleblue] (5.55,4.71) rectangle (7.10,5.32);
\draw[yellow(ncs),line width=2pt               ] (7.55,4.71) rectangle (9.10,5.32);

\end{tikzpicture}\label{fig:bleck3d}}
\caption{Stencil for $z$-stretch example.  Identical stencil values are highlighted with identical boxes
or ovals (dotted rectangles, solid rectangles, dotted ovals, solid ovals)
}
\end{figure}
depicts the matrix stencil for an interior grid point as a function of $\alpha$ for our example $z$-stretched mesh
(see \autoref{sec:stencil3d} for details).
When $\alpha \gg 1$,
one would hope that off-diagonals corresponding to connections within the middle $xy$-plane would
be retained \eqref{eq:crit_sa} and that connections between the central mid-plane vertex and a vertex
in a different $xy$-plane would not.
Figure~\ref{fig:3dsa}
\begin{figure}[htb!]\label{fig:3dsa}
\centering
\includegraphics[trim=2 8 2 22,clip,scale=0.5]{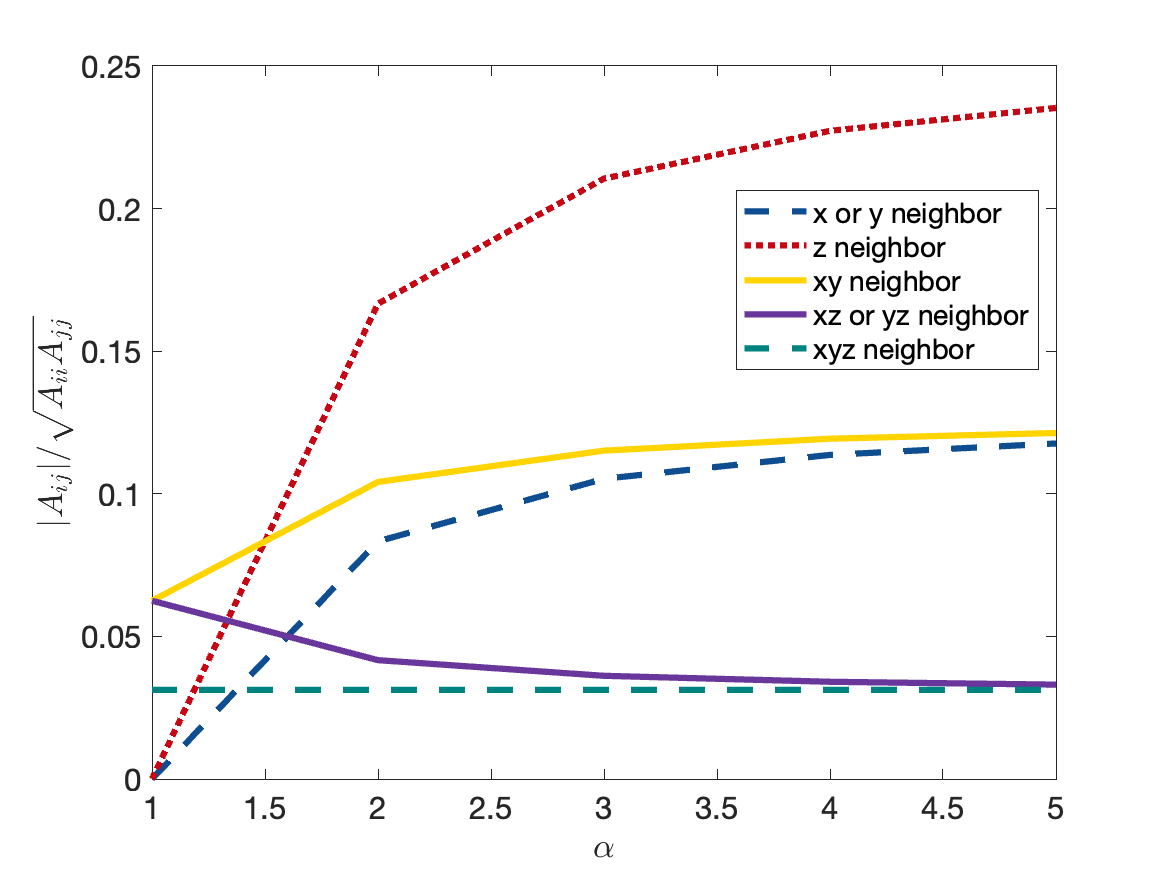}
\caption{Traditional SA criterion values for Figure~\ref{fig:bleez} stencil.}
\vskip -.1in
\end{figure}
plots the scaled matrix entries $|A_{ij}| / \sqrt{A_{ii} A_{jj}}$ which must be less than
the user-supplied $\theta$ in order to be dropped.
For any $\theta$ less than approximately $.1 $, the criterion will properly retain all of the mid-plane entries (those associated
with the blue and yellow lines)
for $\alpha \ge 3$. Additionally,
almost all entries (for the purple and green lines) associated with the top and bottom planes
will be discarded for $\alpha \ge 3$ and $\theta > .037$.
Unfortunately, however, the dotted oval
stencil values (just above
and below the central value) give rise to the largest magnitude scaled entries. Thus, the standard
smoothed aggregation criterion fails completely on this example as it would not restrict coarsening
to just within each $xy$-plane even when $\alpha$ is large.
While these failure could be circumvented by always dropping all positive
off-diagonal entries, this trick would not necessarily work for other PDEs and other discretizations.
In general, we find that basing strong/weak decisions solely on the scaled matrix entries can be unreliable
for finite element discretization due in part to a smearing effect related to the inner products between basis functions.

Though somewhat less obvious from the above plot, there are also shortcomings to how matrix entries
are scaled with the standard smoothed aggregation criterion. Notice that any $\theta$ greater than
$ .0625 $ will discard all off-diagonal entries for $\alpha = 1$ even though relaxation is effective
at removing oscillatory error in all directions. This implies that
we would need to choose $ .037 \lessapprox  \theta \lessapprox .0625 $
in order to successfully use the same $\theta$ over a range of different $\alpha$.
This has implications for finite element applications where the element aspect ratio varies significantly
across the mesh, which is quite common in many cases such as when boundary layers are present.
Such a limited range of suitable thresholds for such a simple model problem
hints at a more problematic concern for more complex applications where users
have little intuition about how best to choose this $\theta$. While part of this shortcoming
is tied to the  misleading nature of the matrix coefficients in this case, another part
of the difficulty arises from how strength-of-connection entries are scaled. It is worth noting that the
classical AMG criterion given by \eqref{eq: tradclass} always retains at least one negative off-diagonal
entry (if present), and so it avoids the pitfall of dropping all off-diagonals entries, which
can easily happen with the smoothed aggregation criterion. One highly problematic case
for the smoothed aggregation criterion is given in \autoref{sec:stencil2d}
for a 2D version of this same stretching problem where
$\theta$ must be less than $.125$ to avoid dropping all off-diagonals in the isotropic case
but needs to be greater than $.125 $ to properly semi-coarsen in the anisotropic
case\footnote{If the trick of discarding all positive off-diagonals is not used, then $\theta$ must be greater than $.25$.}.
As will be discussed in \autoref{sec: dropping}, this issue where all off-diagonals are dropped tends
to be more likely when the number of nonzeros per row is large (which is common for AMG on coarse grids).

There is one additional smoothed aggregation sub-step that affects behavior for stretched-mesh/anisotropic
applications related to the characteristics of the dropped matrix, $\tilde{A}$. 
In smoothed aggregation, a prolongator is constructed in two steps.
The first step builds a tentative prolongator $P_t$ which is effectively a boolean
matrix (values of zero or one) for scalar PDE applications. Specifically,
$$
(P_t)_{ij} = \left \{  \begin{array}{cc} 1 & \mbox{if}~ i \in {\cal A}_j \\
0 & \mbox{otherwise}   \end{array} \right .
$$
where again the set ${\cal A}_j$ includes the index of all fine level vertices in the $j^{th}$ aggregate
and ${\cal A}_j \cap {\cal A}_k = \emptyset $ for all $j \ne k$.
\REMOVE{ More sophisticated tentative prolongators
are used for example with linear elasticity applications in conjunction with a QR algorithm, which we
do not consider in this paper.  It should be noted that a QR algorithm is sometimes applied to post-process
$P_t$'s and improve its condition number, which is motivated by issues for linear elasticity.
However, we find that for scalar PDEs the QR algorithm should be avoided as it can occasionally have a negative impact
on the dropping algorithms on subsequent coarse grids (see Section~\ref{sec: xyz}.
}
As the tentative prolongator is not generally optimal in a standard multigrid V cycle\footnote{
Optimal convergence can be achieved with $P_t$ using special multigrid cycles and small aggregates~\cite{https://doi.org/10.1002/nla.542}.},
a Jacobi-like prolongator smoothing step is then applied to improve the interpolation operator via
\begin{equation} \label{eq: prolong smooth}
P = ( I - {\omega} \widetilde{D}^{-1} \tilde{A} ) P_t
\end{equation}
where $\widetilde{D}$ is the diagonal of $\tilde{A}$ and ${\omega} $ is a damping parameter
that is typically chosen as an inexpensive approximation to $ 4 / ( 3 \rho( \widetilde{D}^{-1} \tilde{A})) $
with $\rho(K)$ defined as the spectral radius of $K$.
The dropped matrix must be used in \eqref{eq: prolong smooth} so that the number of nonzeros
in $P$ is not large when aggregation is applied to
an $\tilde{A}$ where many dropped entries have been excluded from $A$.
In particular, using $A$ in~\eqref{eq: prolong smooth} could
lead to excessive fill-in when constructing the coarse version of the PDE operator
via the Galerkin projection
$P^T A P$.
Additionally,
the support of each of
the prolongator basis function (nonzero pattern within each column) would extend in weak directions
if $A$ is instead used in \eqref{eq: prolong smooth}, and so interpolation would combine information
in oscillatory directions, which we had hoped to avoid.

To maintain smoothed aggregation's optimal convergence properties, it is important that the action
of $\tilde{A}$ applied to a smooth vector should approximately match the action of $A$ applied to
the same smoothed vector, which may not be the case when $\tilde{A}$ is defined by only
zeroing out some entries from $A$. Instead, the retained nonzeros must also be modified
to at least ensure that
\begin{equation} \label{eq: constant action}
s \approx \tilde{s} ~~\mbox{where}~~ s = A v ,~~\tilde{s} = \tilde{A} v, ~~\mbox{and}~~ v
~\mbox{is a constant vector.}
\end{equation}
It is worth
noting that most entries of $s $ are typically zero for many discretized PDE matrices,
mirroring the fact that  differentiation of a constant function
is identically zero.  To guarantee~\eqref{eq: constant action}, the
standard smoothed aggregation procedure \cite{VaMaBr96}
defines $\tilde{A}$ in the following way
\begin{equation}\label{eq: standard lumping}
\tilde{A}_{ij} = \left \{  \begin{array}{cl} A_{ij} & \mbox{if}~ i \ne j ~~\mbox{and}~~ (i,j) \in {\cal S} \\[2pt]
                                          0    & \mbox{if}~ (i,j) \notin {\cal S} \\[2pt]
                                       A_{ii} + \sum_{(i,k) \notin {\cal S}} A_{ij} & \mbox{if}~ i = j
\end{array} \right .
\end{equation}
where ${\cal S} $ is the sparsity pattern that includes all nonzeros of $A$ that are not dropped.
That is, the sum of the dropped entries within each row is lumped to the diagonal entry of that row.

While the standard lumping procedure often works fine, we have noticed scenarios where
for some $i$, $\tilde{A}_{ii} < 0$ even though $A_{ii} > 0$.
Even a single negative diagonal entry can generate
a $P$ (via \eqref{eq: prolong smooth}) that is equivalent to or worse than $P_t$.
Additionally, $\omega$ might be nonsensical (e.g., negative ) due to a negative
eigenvalue estimate that again ruins the prolongator smoothing step.
Though somewhat contrived, it is possible to illustrate a bad lumping scenario
even with our simple stretched example. Specifically,
if the dropped entries from the stencil in~ Figure~\ref{fig:bleez} correspond to those
with values
$ -\frac{1}{2} - \alpha^2$ (green dotted rectangles),
$ -2 - \alpha^2 $ (purple ovals),
and
$1 - 4 \alpha^2$ (yellow rectangles),
then it is easy to verify that $\tilde{A}_{ii} = 0$ using
\eqref{eq: standard lumping}.
\REMOVE {
Notice that with the standard
$SA$ criterion, the green dotted rectangles and the purple ovals
would be dropped for large $\alpha$ over a wide range of thresholds
while for very large $\alpha$ values the scaled entries associated
with blue dotted rectangles and yellow rectangles are quite close.
}
Recalling our earlier remarks to
``{\it avoid coarsening in directions associated with neighboring grid points that are relatively far,"}
we note that in terms of Euclidean distance from the central
point, the vertices associated with the $ -\frac{1}{2} - \alpha^2$ terms
are the furthest for $\alpha > 1$. The vertices for the  $ -2 - \alpha^2 $ terms are
the second furthest, and the vertices for the $1 - 4 \alpha^2$ terms are the third furthest
when $ 1 <  \alpha < \sqrt{2}$.
Though this failure range is somewhat narrow, it illustrates
that a criterion based on distance in conjunction with standard smoothed
aggregation lumping could easily lead to problematic diagonal entries for a subset of $\tilde{A}$'s rows.
A more realistic example will be shown later in the paper.

\section{Smoothed Aggregation Enhancements}\label{sec:new algorithms}
We consider alternatives to the standard smoothed aggregation process.
The first employs a distance Laplacian strength-of-connection (SOC) matrix. The second
considers trade-offs between symmetric and non-symmetric scaling algorithms such as~\eqref{eq: tradclass}. The third alternative classifies
strong and weak connections based on the first {\it large} gap between scaled SOC entries
within a row.  The final enhancement defines $\tilde{A}$ via a
new lumping method.

\subsection{SOC matrix}
The strength-of-connection matrix typically coincides with the discretization matrix $A$. However,
a different auxiliary matrix could have more desirable properties.
In addition to the dropping phase, this auxiliary matrix
could be used within the prolongator smoothing phase. This might even be necessary for saddle point PDEs
if some diagonal entries are identically zero (c.f.,
\cite{doi:10.1137/23M1584514}). 
In this case,
$A$ is
only used in the Galerkin projection step and not used to construct $P$. In this paper, we instead focus
on the scenario where dropping decisions are based on the distance Laplacian while
a dropped version of $A$ is use within the prolongator smoother step (see Section~\ref{sec: sa}), and so $P$
is influenced by $A$'s characteristics.

The {\it distance Laplacian} described in \cite{ddproc06} is
defined as
\begin{equation}\label{eq:entry_dlap}
L_{ij} = -\|\mathbf{x}_i - \mathbf{x}_j\|^{-2}, ~~\mbox{for}~~ A_{ij} \ne 0 ~~\mbox{and}~~ i \ne j
\end{equation}
and
$L_{ii} = -\sum_{i\not=j}L_{ij}$.
Clearly, $L$ has the following characteristics:
\begin{itemize}
\item large (small) magnitude off-diagonal entries correspond to relatively near (far) neighboring mesh vertices from the vertex associated with the diagonal entry;
\item same sparsity pattern as $A$, and $L$ is symmetric if the pattern is symmetric;
\item negative off-diagonal entries;
\item all row sums are zero.
\end{itemize}
For the model problem, SA coarsening applied to the distance Laplacian should mimic the desired semi-coarsening
behavior because off-diagonal entries are the square of inverse
distances.  Notice that distance Laplacian entries resemble those of a 5-pt (or 7-pt) finite difference Poisson discretization (FD) in
two (or three) dimensions on a structured mesh with uniform spacing. In particular, the
off-diagonal entries 
are
$-h_x^{-2}$, $-h_y^{-2}$, and $-h_z^{-2}$
in the $x$, $y$, and $z$ directions respectively, where $h_k$ is the mesh
spacing along the $k^{th}$ axis. In this way, $L$ can be viewed as a very crude approximation
to a Poisson operator that avoids finite element (FE) smearing effects. 
In terms of FE solution accuracy, these smearing effects are addressed by the mass matrix and linear system right hand side, both
not normally available to a typical AMG setup procedure.
\REMOVE{
of a finite discretization operator
in the $x$, $y$, and $z$ directions are $O(h_x^{-2})$, $O(h_y^{-2})$, and $O(h_z^{-2})$ respectively where $h_k$ is the mesh
spacing in the $k$ direction. Thus, the off-diagonal entries with finite differences are effectively the square of the reciprocal distances
while the distance Laplacian entries are just the reciprocal distances. While this might seem like a significant difference,
we note that if the classical AMG dropping criterion~\eqref{eq: tradclass} is used for the distance Laplacian and
for the 5-pt 2D (or 7pt 3D) FD stencil, then the dropping decisions are identical if the FD threshold is the square of the distance Laplacian
threshold. For example, when the mesh spacing in the x-direction is smallest, then the largest negative off-diagonal entry is
$ 1/h_x^2 $ for the FD stencil while it is $1/h_x$ for $L$ while the $z$-direction stencils are $1/h_z^2$ and $1/h_z$, respectively.
This implies that the dropping decision for the $z$-direction off-diagonal is based on $1/h_z^2 \ge \theta^2 1/h_x^2 $ in the FD case while it is
based on $1/h_z \ge \theta 1/h_x $ in the distance Laplacian case.
}
\REMOVE{
Adapting the standard SA criterion \eqref{eq:crit_sa} to the distance Laplacian
gives the dropping criterion
\begin{equation}\label{eq:crit_dlap}
|L_{ij}| \geq \theta \sqrt{L_{ii} L_{jj}},
\end{equation}
which then defines the sparsity pattern of $\tilde{A}$.
}

One disadvantage of the distance Laplacian
is that coordinates must be supplied by the application.
In our experience, we find that providing
coordinates is easy for most
application developers.
This is in contrast with supplying
complete mesh information (e.g., element geometry/connectivity), which is onerous
and requires non-trivial code interfaces.
Within the AMG implementation, $L$ can be explicitly formed and projected
to coarse levels. Storage can be reduced using a low precision $L$
and/or retaining only its sparsity pattern once dropping decisions are made.
Alternatively, one can implicitly define $L$ by projecting
coordinates to coarser levels, reducing storage further.

A second disadvantage to the distance Laplacian
is that important PDE features might be missed when
coarsening or when defining $P$'s sparsity pattern.
For example, if the PDE is anisotropic or if its
coefficients include large variations or discontinuities
(e.g., due to multiple materials), these features will not
influence coarsening. However, if sub-par AMG performance is mainly due
to mesh stretching for a particular application,
then distance Laplacians
are attractive. For some 
engineering
applications, 
distorted meshes with poor element aspect ratios
are the primary source
of AMG difficulties. It is important to note that
the prolongator 
still incorporates
information from the discretization matrix
through $\tilde{A}$, 
albeit not for coarsening decisions.

\subsection{SOC Scaling and Nonsymmetry} \label{sec: dropping}
SOC scaling is generally necessary. This can be seen by
considering a linear system $ B u = f_b $ where
$ B = \gamma A $ and $f_b = \gamma f$.
This might occur if the mesh spacing is stretched uniformly in all
coordinate directions.
We would expect that dropping
decisions for $B$ would be identical to those when solving
$ A u = f$, but scaling must be introduced to
preserve this invariance. 
\REMOVE{
Clearly, this small modification shouldn't affect
dropping decisions nor the choice of $\theta$, but this will require
scaling the SOC matrix to
}
Scaling effectively defines diagonal matrices ${D}_\ell$ and ${D}_r$ and considers a scaled
SOC matrix
\begin{equation} \label{eq: generic scaling}
{D}_\ell^{-1} {S} {D}_r^{-1} .
\end{equation}
\REMOVE{Notice that the scaling associated with the
standard SA criterion and the classical AMG criterion are both
of the form ~\eqref{eq: generic scaling}. Specifically,
}
For SA,  ${D}_\ell  = {D}_r$ and $({D_r})_{ii} = \sqrt{S_{ii}}$.
For classical AMG, ${D}_r = I$ and $({D}_\ell)_{ii}
= \max_{k\not=i} -S_{ik}$.
In either case,
$ {D}_\ell$ and ${D}_r $ offset the uniform scaling of  $S$ by $\gamma$.

\begin{figure}[h!]
\begin{center}
\includegraphics[trim = 0.6in 1.1in 4.0in 5.2in, clip,scale=0.5,angle=0,origin=c]{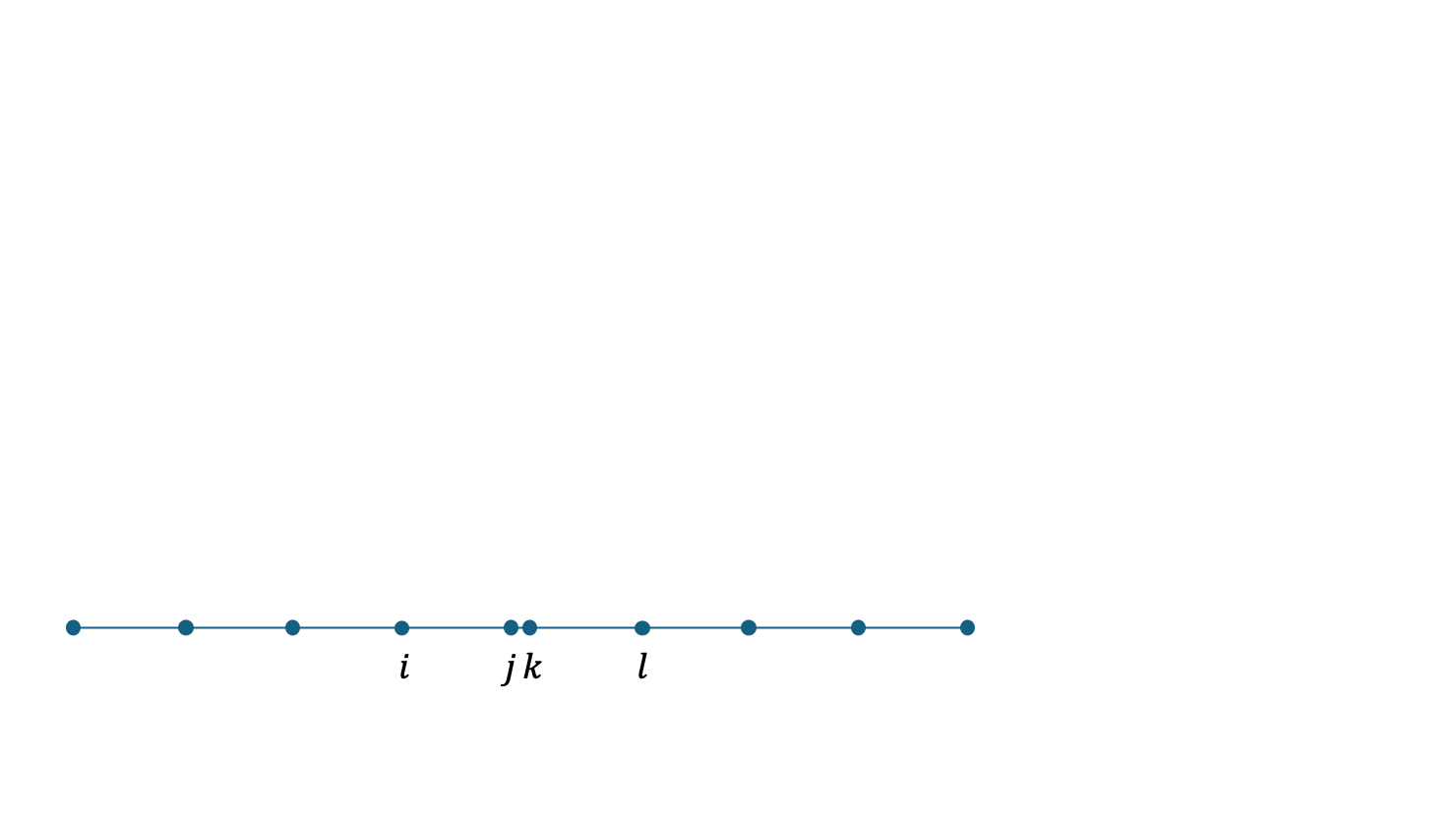}
\vspace*{-.05in}
\caption{Mesh with 2 nearby vertices relative to distances between
other adjacent vertices.                          } \label{fig: 1D stretched}
\end{center}
\end{figure}

Our original preference for the smoothed aggregation criterion has been based on symmetry advantages.
Specifically, the graph of the scaled strength-of-connection matrix is represented
with undirected edges when it is a symmetric.  Recalling the aggregation process, any
two vertices in the $m^{th}$ aggregate are connected by a path in the strong scaled strength-of-connection
graph.  When ${D}_\ell^{-1} {S} {D}_r^{-1}$ is non-symmetric, directed edges are needed. Two vertices in ${\cal A}_m$
should be connected by directed edges in the strong graph, but the path does not necessarily respect strong edge orientations.
In fact, the specifics of the aggregate might depend on implementation details.
\REMOVE{
of edges connecting $i$ and $j$ depends
More recently, however, we find
that non-symmetric approaches 
offer important benefits. 
Before introducing the {\it cut-drop} approach,
we briefly discuss the role that symmetry plays in aggregation.
to represent it. In this case, two vertices $i$ and $j$ in ${\cal A}_m$ must be connected by
directed edges, but
there is not necessarily 
a path from $i$ to $j$ that respects edge directions.
}
To understand this, consider the 1D mesh in Figure~\ref{fig: 1D stretched} where vertices
are uniformly spaced except for two vertices that are close together. When $S = L$,
the SA criterion determines that edge $(j,k)$  is strong while
edges $(i,j)$ and $(k,l)$ are weak if a suitable $\theta$ is chosen.
This means that an aggregate
that includes only $j$ and $k$ is always formed.
With the
classical AMG criterion, however, the directed edge $(i \hskip -.04in \rightarrow \hskip -.04in j)$ is strong
as $i$'s left and right neighbors are equi-distant form $i$'s location. Conversely,
edge $(j \hskip -.04in \rightarrow \hskip -.04in i)$ is weak as $j$'s right neighbor is much closer than its left neighbor. The same sort of split decision occurs
for $(k \hskip -.04in \rightarrow \hskip -.04in l)$ and $(l \hskip -.04in \rightarrow \hskip -.04in k)$. This implies that vertices $j$ and $k$
will reside in different aggregates if the aggregation procedure initially chooses the $i^{th}$ vertex
to form one new aggregate (which includes $j$ and $i$'s left neighbor)
and the $l^{th}$ vertex to form another (which includes $k$ and $l$'s right neighbor). Notice that the path from
$j$ to $i$'s left neighbor relies on strong edge $(i \hskip -.04in \rightarrow \hskip -.04in j)$, which has the wrong orientation.
If instead vertex $j$ or $k$ happens to be chosen as the initial vertex to form an aggregate, then an aggregate that includes only
$j$ and $k$ is formed as in the SA criterion case.
In this sense, a symmetric scaled SOC matrix seems superior
in that it always aggregates the nearby $j$ and $k$ vertices together regardless of aggregation details.
\REMOVE {
While abrupt jumps in mesh spacings are not generally advised, we note that
applications with discontinuities can lead to similar situations (e.g., an
application where a thin wire is surrounded by air). When addressing a thin wire,
it is best to avoid aggregating its vertices with adjacent non-wire vertices, which smears the wire on coarse meshes.
}

So what are the disadvantages of the symmetric SA criterion? The biggest disadvantage is that it can frequently determine that
all off-diagonal nonzeros should be dropped even when relaxation is effective in all directions. To examine
this, consider the case where all of $S$'s off-diagonal entries are negative (e.g., using
many common finite difference schemes or using a distance Laplacian as the SOC matrix).
The sum of entries within a row
corresponding to the domain interior is frequently zero.
In this case, the SA strong criterion can be re-written as
\begin{equation} \label{eq: average}
\frac{|S_{ij}|}{\sqrt{p_i p_j} \sqrt{ \left | \bar{S}_{i} \bar{S}_{j} \right | } } \ge \theta
\end{equation}
where $\bar{S}_{k}$ is the average of all off-diagonal nonzeros in the $k^{th}$ row
and $p_k$ is the number of off-diagonal nonzeros in the $k^{th}$ row and column.
When the variance of the off-diagonal values is relatively small (e.g., element aspect ratios close to one),
$S_{ij}/\bar{S}_i \approx 1$ and $S_{ij}/\bar{S}_j \approx 1$.
This implies that the magnitude of the scaled entries is sensitive to the number of nonzeros
within a row/column and that the left hand side of \eqref{eq: average} becomes small for all
off-diagonals when $\sqrt{p_i p_j}$ is large.

In general, the number of nonzeros per row is higher for 3D domains
than for 2D domains.  Even on the finest grid,
the number of nonzeros per row can also vary noticeably depending on how many elements
share a corner vertex.  Further, the number of nonzeros per row is often significantly higher
and more variable on coarser AMG grids.  All of this suggests
that a row might contain only weak off-diagonals simply because the number of nonzeros
in this row is large.  While
software checks can be included to prevent division by zero (as the associated diagonal entry of $\tilde{A}$ must be set to zero to maintain row sums), multigrid convergence and the cost per iteration
(due to a lack of coarsening) can still be seriously impacted when this {\it all-dropping}
occurs across many rows.
Notice that all-dropping cannot happen with the classical criterion as it is guaranteed to
identify at least one off-diagonal in each row
as strong (if there is at least one negative entry). Further, it will properly retain most of the entries in the isotropic
or low variance case. Thus, the non-symmetric criterion has an important advantage over the SA criterion. 

Before concluding this sub-section,
we briefly note that any non-symmetric criterion can be modified
to form a symmetric criterion. One simple strategy classifies an edge $(i,j)$ as strong only if the non-symmetric criterion
classified either $(i \hskip -.04in \rightarrow \hskip -.04in j)$ or $(j \hskip -.04in \rightarrow \hskip -.04in i)$ as strong.
Alternatively, one could only classify edge $(i,j)$ as weak only if the non-symmetric criterion
classified either $(i \hskip -.04in \rightarrow \hskip -.04in j)$ or $(j \hskip -.04in \rightarrow \hskip -.04in i)$ as weak.
Similarly, one could average the scaled entries associated with $(i \hskip -.04in \rightarrow \hskip -.04in j)$ and $(j \hskip -.04in \rightarrow \hskip -.04in i)$ before classification.
\REMOVE{
Unfortunately, the first (second) approach would (might) classify all neighbor connections
as strong in our Figure~\ref{fig: 1D stretched} example.
We generally
find that it is often best to instead  classify $(i,j)$ as weak if the non-symmetric criterion
classified either $(i \hskip -.04in \rightarrow \hskip -.04in j)$ or $(j \hskip -.04in \rightarrow \hskip -.04in i)$ as weak.
That is, it is better to favor semi-coarsening to assure rapid convergence rates (at the expense of iteration cost)
when there is a split decision. 
This would properly treat the Figure~\ref{fig: 1D stretched} case.
However, this
no longer guarantees that at least one off-diagonal within a row is retained.
We consider one symmetrization in this manuscript for one example case.
}
We do not considered symmetrization further in this manuscript as we have not observed a significant advantage to
it on stretched mesh Poisson problems.
\REMOVE{

The scaled off-diagonal entries would be either
$|A_{ij}|/\sqrt{A_{ii} A_{jj}}$
or
$|L_{ij}|/\sqrt{L_{ii} L_{jj}}$
depending on whether or not the SOC matrix is chosen as $A$  or $L$.
Recalling that all off-diagonal entries of $L$ are negative and that
all rows have a sum of zero, we can rewrite the criterion
for the distance Laplacian as
where we have additionally assumed that all rows have the same number of off-diagonal nonzeros
equal to $m$ and
row.

This could again cause all off-diagonal entries to be dropped as discussed in Section~\ref{sec: sa}.

While symmetric scaling is nice in that $\tilde{A}$ maintains the symmetry of $A$ when it is symmetric,
it is worth noting that the Jacobi

it is clear that as the

As noted, the standard criterion can sometimes drop all off-diagonal entries
even in cases when relaxation is effective at smoothing error in all directions.
While we emphasized this point when the SOC matrix is simply taken as the $A$ matrix,
it can also occur when $L$ is used as the SOC matrix especially when the number
of nonzeros per row is large.
}
\subsection{Classification Based on Gaps}
The scaled SOC entries within a row are often distributed in a clustered fashion with several
similar-magnitude values followed by a gap and then additional similar-magnitude values.
For example, entries associated with the normal and tangential directions of a boundary layer often
lie in distinct clusters. A threshold approach must employ a suitable $\theta$ to effectively separate
clusters. In some cases, it can be challenging to find a $\theta$ that does this over
an entire mesh's rows. The cut-drop algorithm (summarized in Algorithm~\ref{alg:cut})
that we now present takes a different approach
that directly aims to find a gap between clusters.
While we have not described or systematically evaluated it previously,
the cut-drop algorithm
was incorporated in the MueLu \cite{MueLu} software several years ago.
\REMOVE{
The cut-drop algorithm is a non-symmetric criterion, but it does include
some symmetry aspects not present with the classical AMG criterion
that allow it to properly address our Figure~\ref{fig: 1D stretched} example.
Like the classical criterion it guarantees that at least one
off-diagonal is retained (barring rows with no off-diagonals). It differs in one further fundamental way
from the classical criterion in that classification decisions
revolve around finding a gap in the values of the scaled SOC
entries within a row.
}
\begin{algorithm}[htb]
\begin{algorithmic}
\STATE $E_r \leftarrow \emptyset$
\STATE \textcolor{red}{/* Define symmetric-scaled values */}
\STATE $v_{ij} \leftarrow |S_{ij}|/\sqrt{S_{ii} S_{jj}}~~~\forall i\ne j ~\mbox{and}~ (i,j) \in {\cal G}_s $ \textcolor{red}{/* ${\cal G}_s $ is graph of $S$ */}
\STATE
\vskip -.1in
\FOR{$i = 1,\ldots,n$}
  \STATE $(\hat{v},\hat{j}) \leftarrow $\texttt{descendSort}$(v_{i*})$ \textcolor{red}{/* $\hat{j}_k$ is column index of $\hat{v}_k$ ($k^{th}$ largest nonzero).*/}
  \STATE $E_r \leftarrow E_r \cup (i,\hat{j}_1) $
  \STATE $ k  \leftarrow 2 $
  \STATE  \texttt{lastEdgeStrong} $\leftarrow true $

  \WHILE{( $ (k \le $\texttt{length}$(\hat{v})$ AND \texttt{lastEdgeStrong})}
     \IF {$\hat{v}_{k} / \hat{v}_{k-1} \geq \theta_g$}
       \STATE $E_r \leftarrow E_r \cup (i,\hat{j}_k)$
     \ELSE
       \STATE  \texttt{lastEdgeStrong} $\leftarrow false$
     \ENDIF
     \STATE $ k  \leftarrow k + 1 $
  \ENDWHILE
\ENDFOR
\end{algorithmic}
\caption{${E}_r$ = CutBasedDropping($S,\theta_g$)\\
\textcolor{red}{/* $S \in  \mathbb{R}^{n} \times \mathbb{R}^{n}$: SOC matrix;~~}
\textcolor{red}{$\theta_g$: gap tolerance;}
\textcolor{red}{~~${E}_r$: retained edges \hskip .1in */} \label{alg:cut}}
\end{algorithm}
Cut-drop begins by scaling entries symmetrically (like standard SA).
Instead of thresholding the scaled $v_{ij}$ entries,  
they are sorted within each row. 
The largest scaled entry within a row is always classified as strong (like the classical criterion).
The next largest $v_{ij}$ within a row is then
classified as strong only if the ratio between this $v_{ij}$
and the previous one is above (or equal to) a tolerance. If this
is the case, the process is repeated for the next entry computing
a ratio with the most recently classified strong entry. Once
an edge is classified as weak, then all remaining smaller scaled entries are
classified as weak within the row. Notice that while the scaling is symmetric, classification
within each row is independent of other rows, and so the final strong matrix is non-symmetric.

\subsection{Lumping the dropped matrix}
To maintain matrix row sums while also preserving the sign of the diagonal entry, a new lumping procedure is proposed
and summarized in Algorithm~\ref{alg:lump}.
\begin{algorithm}[htb]
\begin{algorithmic}
\STATE \textcolor{red}{/* Construct initial dropped matrix $\bar{A}$ */}
\STATE $\bar{a}_{ij} \leftarrow a_{ij} ,~~ \forall (i,j) \in \tilde{\cal S} $
\STATE \texttt{rowSum} $\leftarrow A \vec{1} $ \textcolor{red}{/* $\vec{1}$ is vector with all entries equal to one. */}
\STATE $\delta$\texttt{RowSum} $\leftarrow $~\texttt{rowSum} $ - \bar{A} \vec{1} $
\STATE $\bar{A}$\texttt{RowAbsSum} $\leftarrow |\bar{A}| v $
\STATE \textcolor{red}{/* Update dropped matrix $\tilde{A}$ */}
\STATE $\tilde{a}_{ii} \leftarrow \bar{a}_{ii} + ~\delta$\texttt{RowSum}$_i                                                 ,~~ \forall (i,i) \in \tilde{\cal S} ~|~
\delta$\texttt{RowSum}$_i > 0$
\STATE $\tilde{a}_{ij} \leftarrow \bar{a}_{ij}
+~ \delta$\texttt{RowSum}$_i
~(|\bar{a}_{ij}|/\bar{A}$\texttt{RowAbsSum}$_i) ,~~ \forall (i,j) \in
\tilde{\cal S} ~|~ i\not=j ~\cap~ \delta$\texttt{RowSum}$_i < 0$
\end{algorithmic}
\caption{$\tilde{A}$ = DropWithLumping($A, \tilde{\cal S}$)\\
\textcolor{red}{/* $A$: undropped matrix;~~ }
\textcolor{red}{   $\tilde{\cal S}$: target sparsity pattern;~~}
\textcolor{red}{   $\tilde{A}$: dropped matrix } \label{alg:lump}}
\end{algorithm}
The new algorithm is much simpler than the procedure provided in~\cite{savariants22}
that has the
same goals but lacked guarantees concerning the sign of the diagonal entries of the resulting dropped matrix.
In addition to this guarantee, 
we find that the new
algorithm works better. The idea is to perform traditional diagonal lumping when it increases the magnitude
of the diagonal entry. If traditional diagonal lumping would instead lower the diagonal entry, we
instead distribute the lumped quantity over
all the retained entries in a manner that is proportional to their magnitude. Recalling that
the actual entries of $\tilde{A}$ are only used during the prolongator smoothing step and that
within this step $\tilde{A}$'s entries are left scaled by $\widetilde{D}^{-1}$,
we are interested in two properties  related to the lumping procedure.
The first concerns the sign of the lumped diagonal entry
while the second centers on the magnitude of the negative and of the positive entries in the matrix
$\widetilde{D}^{-1} \tilde{A}$.

\begin{theorem} \label{thm: positive diagonal}
Let $A$ be an $n \times n$ matrix with
sparsity pattern ${\cal S}$, having positive diagonal entries and
and non-negative row sums.
Consider a matrix $\tilde{A}$ whose sparsity pattern
$\tilde{\cal S} \in {\cal S}$ whose nonzero entries (i.e., $(i,j) \in \tilde{\cal S}$) are given by
\begin{equation} \label{eq: Atilde modification}
\tilde{a}_{ij} =
\left \{
\begin{array}{ll}
a_{ij} + e_i \frac{ |a_{ij}|}{\sum_{k\in {\cal \tilde{S}}} |a_{ik}|}  \hskip .2in &\mbox{for}~~ e_i < 0 \\[2pt]
a_{ij}                                                                    &\mbox{for}~~ e_i \ge 0  ~\mbox{and}~ i\ne j\\[2pt]
a_{ij} + e_i                                                              &\mbox{for}~~ e_i \ge 0  ~\mbox{and}~ i = j\\
\end{array}
\right .
\end{equation}
where
\begin{equation}
e_i = \sum_{k ~|~ (i,k)\in {\cal S}\setminus\tilde{\cal S} } a_{ik}.
\end{equation}
That is, $e_i$ is the sum of nonzero entries in $A$ excluded from the $i^{th}$ row of
$\tilde{A}$.
Additionally, assume that for any row $i$, there exists at least one $(i,k) \in
\tilde{\cal S}$ such that $k\not=i$ and $a_{ik}<0$ (at least one negative off-diagonal entry in each
row of $A$ is retained in $\tilde{A})$.
Then, the following two properties hold:
\begin{itemize}
\item For all $(i,j) \in \tilde{\cal S}$, $sgn(a_{ij}) = sgn(\tilde{a}_{ij})$ (Entries retain their sign).
\item For all $i$, $\sum_j a_{ij} = \sum_j \tilde{a}_{ij}$ (Identical
row sums).

\end{itemize}
\end{theorem}
\begin{proof}
We consider first the case where $e_i \ge 0$.  As this approach shifts
positive values to the diagonal, all signs are preserved as is the row sum.
This corresponds to a traditional lumping procedure that does not
alter the retained off-diagonal values and modifies the diagonal to maintain the row sums.

Now we consider the $e_i < 0$ case.  We first observe that we can rewrite the expression for $\tilde{a}_{ij}$ as
\begin{equation} \label{eq: scaled form}
\begin{array}{ll}
\tilde{a}_{ij} = {a}_{ij}\left ( 1
+ \frac{e_i}{\sum_{k \in \tilde{\cal S}_i} |a_{ik}|} \right ) \hskip .2in & \mbox{where}~ a_{ij} > 0 \\[3pt]
\tilde{a}_{ij} = {a}_{ij}\left ( 1
- \frac{e_i}{\sum_{k \in \tilde{\cal S}_i} |a_{ik}|}\right )
& \mbox{where}~ a_{ij} < 0 .
\end{array}
\end{equation}

The parenthesized expressions in \eqref{eq: scaled form}
include only the sums of excluded or retained
entries (or their absolute values), and so these scaling factors are the same
for each entry with that particular sign within a given row $i$.
For  $a_{ij} < 0$, we have $e_i < 0$ and a positive denominator,
making the parenthesized expression positive so the the first property
holds.

For the $a_{ij} > 0$ case, we define
\begin{eqnarray}
r_i^{(p)} &=& \sum_{j~|~ (i,j) \in \tilde{\cal S} ~\mbox{and}~ a_{ij} > 0} a_{ij} ~\mbox{(sum of positive retained entries)},\\
r_i^{(n)} &=& \sum_{j~|~ (i,j) \in \tilde{\cal S} ~\mbox{and}~ a_{ij} < 0} a_{ij} ~\mbox{(sum of negative retained entries)},\\
s_i &=& e_i + r^{(p)}_i  + r^{(n)}_i ~\mbox{(row sum)},\label{eq:proof_si}
\end{eqnarray}
for each row where $e_i < 0$.  We note that $s_i = \sum_j a_{ij}$, the sum of $A$'s entries in the $i^{th}$ row.
We can now rewrite \eqref{eq: scaled form} as
\begin{equation}
\tilde{a}_{ij} = {a}_{ij}\left ( 1 + \frac{e_i                  }{r_i^{(p)} - r_i^{(n)}} \right ),
\end{equation}
for all positive $a_{ij}$, which obviously includes the diagonal
entry.  For the third property to hold, we need
\begin{equation} \label {eq: sum fraction}
1 + \frac{e_i                  }{r_i^{(p)} - r_i^{(n)}} > 0
\end{equation}
or equivalently that
\begin{equation}
e_i  + r_i^{(p)} - r_i^{(n)} > 0\label {eq:sum_fraction2}
\end{equation}
as the denominator in \eqref{eq: sum fraction} is always positive.
Substituting \eqref{eq:proof_si} for $e_i$,
into \eqref{eq:sum_fraction2} the condition becomes
\begin{equation} \label{eq: row sum condition}
 s_i - 2 r_i^{(n)} > 0 .
\end{equation}
We note that $r_i^{(n)} < 0$ because we assumed that at least one negative off-diagonal is retained.
Since $s_i$ is non-negative, it follows that all positive retained
entries (including the diagonal entry) remain positive, and the first
property holds.

For the second property in the $e_i<0$ case, the row sum of the $i^{th}$ row of $\tilde{A}$ is obtained by adding all entries in the top line
of \eqref{eq: Atilde modification} and corresponds to
\begin{equation*}
r_i^{(p)} + r_i^{(n)} + e_i                    \frac{\sum_{k\in \tilde{\cal S}} |a_{ik}|}{\sum_{k\in \tilde{\cal S}} |a_{ik}|}  .
\end{equation*}
As this is just the definition of $s_i$, we have shown
that $\tilde{A}$ has the same row sums as those of $A$, completing the
proof.
\end{proof}
\vskip .1in
\noindent
\begin{remark}
Notice that if $s_i = 0$ and no negative entries are retained, then any lumping strategy that
maintains row sums must lead to a diagonal entry that is either
identically zero or negative. This highlights the importance of
keeping at least one negative off-diagonal value.
Further notice that
should a row have a negative sum (perhaps due to a somewhat atypical boundary condition), then it is still
possible that the diagonal remains positive if the sum of the retained negative entries is large enough.
Should \eqref{eq: row sum condition} not hold (and so the sign of the diagonal entry changes),
we note that traditional lumping would also change the sign of the diagonal as it adds
a larger magnitude negative value to the diagonal than the top expression in \eqref{eq: Atilde modification}.
Finally, we note that additional strategies could be considered
for rows with negative row sums (e.g., lumping only to off-diagonal negative entries), but
we have not experimented with this as we
have not encountered a situation where \eqref{eq: row sum condition} does not hold in our experiments.
\end{remark}

\begin{corollary}\label{cor:a}
For rows where $e_i < 0 $, the positive entries of $      {D}^{-1}       {A}$ are identical to those of
$      \widetilde{D}^{-1} \tilde{A}$ when \eqref{eq: Atilde modification} is used
to define $\tilde{A}$ and the same assumptions in Theorem~\ref{thm: positive diagonal} hold.
Likewise for rows where $e_i < 0 $, the negative entries of $      \widetilde{D}^{-1} \tilde{A}$ are
obtained by scaling the negative retained entries of $      {D}^{-1}       {A}$
by the factor
\begin{equation} \label{eq: negative scaling factor}
\frac{1 - \frac{e_i}{\sum_{k \in \tilde{\cal S}_i} |a_{ik}|} }
     {1 + \frac{e_i}{\sum_{k \in \tilde{\cal S}_i} |a_{ik}|}}
\end{equation}
which is greater than one.
\end{corollary}
\begin{proof}
From the proof of Theorem~\ref{thm: positive diagonal}, we established that all
positive entries within a row are scaled by the same factor
given by the upper expression
between parenthesis of \eqref{eq: scaled form}
and hence
all positive entries of $      {D}^{-1}       {A}$ and
$\widetilde{D}^{-1} \tilde{A}$ must coincide when $e_i < 0$.
We also established in the Theorem~\ref{thm: positive diagonal} proof that
all negative entries within a row are scaled by
the lower expression between parenthesis in \eqref{eq: scaled form}, and
so~\eqref{eq: negative scaling factor} follows.
Additionally, the denominator is positive and must be smaller than the numerator as
$e_i < 0$. Thus, the scaling factor must be greater than one.
\end{proof}

\begin{remark}
While it is not completely clear what properties are best for an
ideal lumping procedure, we find it natural that the only difference between
$      {D}^{-1}       {A}$ and $\widetilde{D}^{-1} \tilde{A}$ for rows where $e_i < 0$ is that
the constant scaling of the negative entries is a bit larger than the
positive entries (which are not scaled) as this larger
scaling can be viewed as a way to compensate
for the fact that more negative information has been dropped
from the stencil, which is reflected by the sum of the excluded entries being negative.
Specifically, it can be easily shown that the scaling factor for negative entries when
$s_i = 0$ (i.e., interior rows) and $e_i < 0$ is given by $ r_i^{(p)} / |r_i^{(n)} | $. That is,
the scaling factor only becomes large to compensate for the case when most negative entries
have been dropped.
\end{remark}

\begin{remark}
While preserving the sign of the diagonal entries, the lumping
procedure described in \autoref{thm: positive diagonal}
does not maintain the symmetry of $A$ nor does it guarantee that $\tilde{A}$ has all positive
eigenvalues (which is also not guaranteed by traditional diagonal lumping). While this can cause
degradation in the quality of the prolongator smoothing step, we find it far less harmful
than when even a single diagonal entry changes sign.
\end{remark}

%


\section{Computational Examples}\label{sec:examples}
A variety of stretched mesh examples are considered to assess
different sub-stage classification methods.  All examples employ an
algebraic multigrid preconditioner through the MueLu \cite{MueLu} software package.
 MueLu is set to generate  multigrid levels until the coarse problem is under
1,000 unknowns. A multigrid V cycle is used with 1 pre- and 1 post-
symmetric Gauss-Seidel iteration on all tests.
For all examples, linear finite elements using either quadrilateral,
triangular (only for the cavity mesh)  or hexahedral elements are
employed to discretize a Poisson operator. When solving the associated linear systems,
the right hand side is chosen (for all but the cavity problem)
so that the 2D solution is $u(x,y) = 1+x+y+xy$,
the 3D solution is $u(x,y) = 1 + x + y + z + x y + x z + y z + x y z$, and the
initial guess satisfies Dirichlet boundary conditions, but is zero elsewhere.

Meshes for the examples described in the first three sub-sections are generated using
Pamgen~\cite{Pamgen}.  Appendices \ref{sec:pamgen2d}-\ref{sec:pamgen3d_radtri}
provide the relevant input decks. Meshes in Section~\ref{sec:results_2d_cavities}
are generated using Cubit~\cite{CUBIT} in conjunction with post-processing
to add boundary layers. Finally, the Section~\ref{sec:application meshes} meshes
were generated by different application teams at Sandia.

\REMOVE {
In Sections~\ref{sec:results_2d}-\ref{sec:results_3d_cone},
We consider solving \eqref{eq:poisson}
using the TrilinosCouplings package of Trilinos
\cite{Trilinos-Overview} to generate discretizations.
The remaining examples use a finite element code implemented in
Matlab to generate matrices. The first of these focuses on boundary
layer meshes with cavities. Another problem uses meshes from
a Sandia hypersonic code. Even though our focus has been on mesh stretching
a final example is given for a
material variation problem  where the mesh
is not stretched by the problem contains discontinuities in the materials.

We conclude with two final examples that employ a Matlab version of
smoothed aggregation on meshes that include both quadrilateral
and triangular elements. While the first of these considers a constant



We conclude with two final examples that employ a Matlab version of
smoothed aggregation on meshes that include both quadrilateral
and triangular elements. While the first of these considers a constant
coefficient Poisson PDE, the last is a material variation
Poisson example.
\RST{Could just dump matrices? What about symmetry?}
}

\subsection{2D Tensor-Product Stretched Poisson}\label{sec:results_2d}

\begin{figure}\label{fig:sample_mesh2d}
\begin{center}
\scalebox{-1}[1]{ \includegraphics[scale=0.3,angle=90]{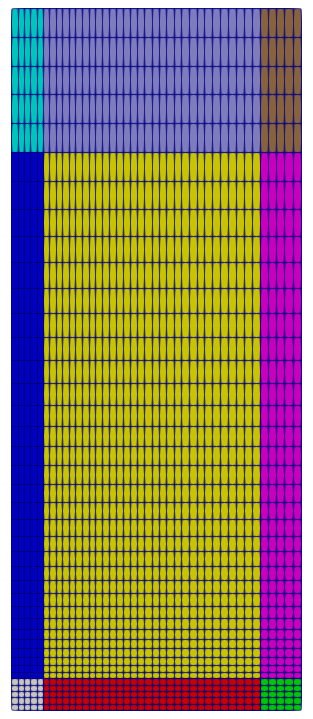}}
\end{center}
\caption{2x de-refinement of 2D tensor-product mesh with stretch factors of $\gamma_1 = 4.5459$ and $\gamma_2 = 1.2877$.}

\end{figure}

A mesh 
with nine regions 
is shown in Figure~\ref{fig:sample_mesh2d}.
It has an isotropic bottom-left region, a uniformly stretched top-right region, and
graded regions in between.
%
Twenty different stretch factors ($\gamma_1$ and $\gamma_2$) logarithmically spaced
between 0.5 and 200 for each dimension are investigated.  
We only consider meshes where the $y$-stretching is greater than or equal to that of the $x$ stretching,
leading to 210 total problems.
Dirichlet conditions are
applied on the bottom boundary while all other boundaries have Neumann conditions.

\begin{figure}\label{fig:live2d}
\hskip -.35in \includegraphics[trim = 2.2in 0.9in 1.9in 0.6in, clip,scale=0.235]{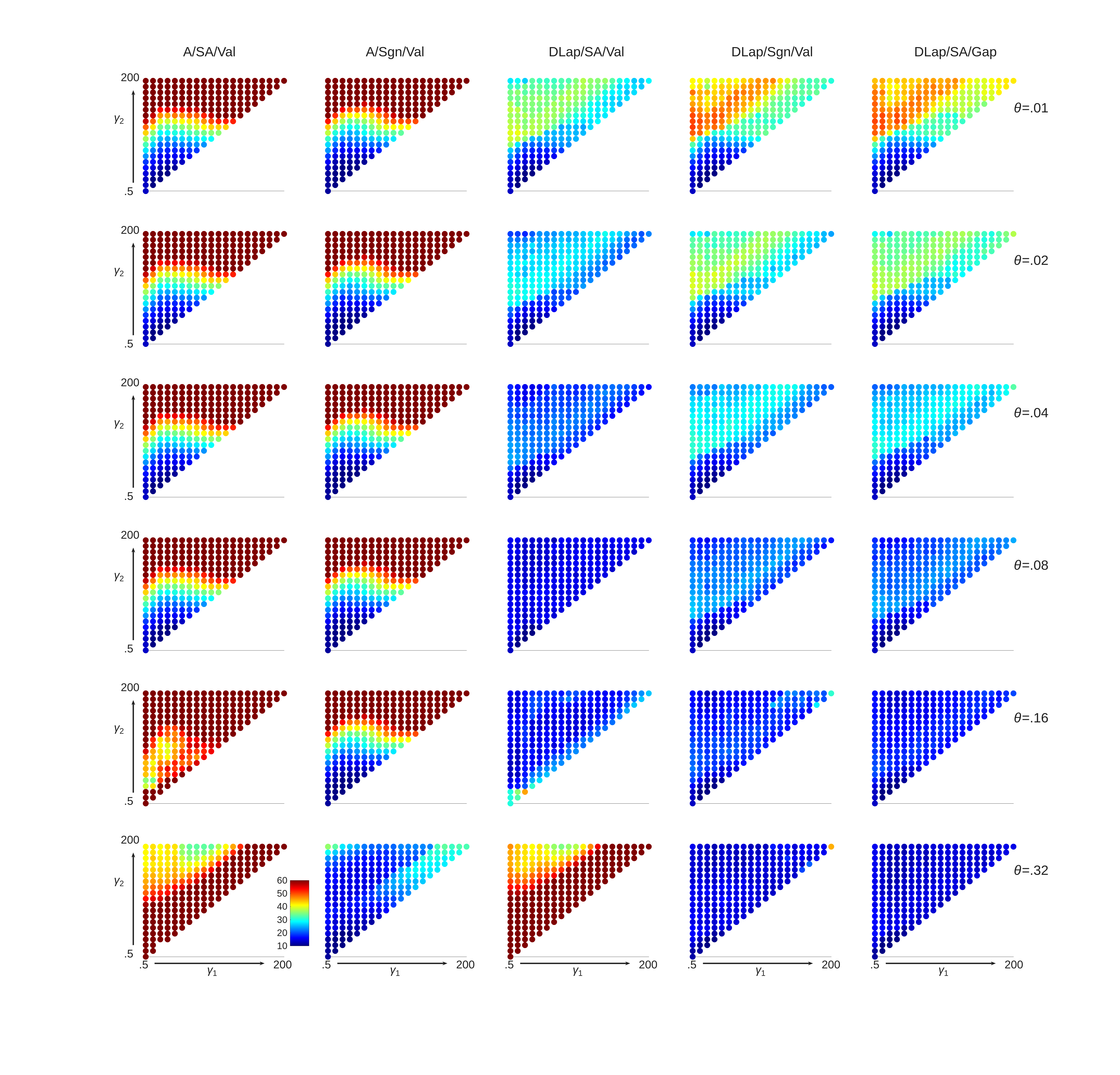}
\vskip -.2in
\caption{Iterations{\hskip .01in}$\times${\hskip .01in}operator{\hskip-.02in} complexity, with the color scale from 10 to 60 for the 2D mesh
Columns represent algorithm choices (described in
Section~\ref{sec:results_2d}) and rows represent tolerances,
$\theta$. Diagonal lumping is use for all experiments.}
\end{figure}

Figure~\ref{fig:live2d} shows results for five algorithm combinations using an AMG-preconditioned conjugate gradient (CG) solver,
stopping when the residual is reduced by 10 orders of magnitude. Each algorithm is denoted by a name of the form
SOC matrix/scaling scheme/classification scheme. The SOC matrix could be either {\sf A} or {\sf DLap}
to indicate use of the matrix $A$ or the distance Laplacian respectively.
The scaling algorithm could be either {\sf SA} or {\sf Sgn} to indicate use of the SA symmetric scaling
scheme or the signed classical AMG scheme respectively.
The classification algorithm could be either {\sf Val} or {\sf Gap} to indicate that classification is based on
scaled SOC entry values or gaps (i.e., cut-drop algorithm). These results all employ a traditional
diagonal lumping scheme.
Each dot corresponds to a single example and measures the overall cost defined as the number of
iterations multiplied by operator complexity. The AMG operator complexity is the ratio
between the sum of the number of nonzeros in the matrix $A$ on all levels divided by the number of $A$'s nonzeros on the finest level.  Blue is the best (cost $= 10$) and red is the worst (cost $= 60$).  To better distinguish between successful runs, we clip colors at $60$. That is, all runs with cost greater than $60$ also appear as red dots.
A missing dot indicates that either AMG did not converge or the setup procedure broke down (e.g., all
off-diagonal entries dropped).  In general, it is clear that the distance Laplacian SOC matrix
is advantageous when compared to using the $A$ matrix as the solver converges well over a wide range of $\theta$'s.
However, {\sf A/Sgn/Val} for $\theta = .32$ also performs well.
Here, the automatic classification of positive off-diagonals as weak plays an important role.
It is interesting to note that {\sf A/Sgn/Val}, {\sf DLap/Sgn/Val}, {\sf Dlap/SA/Gap} all do well with
$\theta = 0.32$. This $\theta$ value is much larger than standard AMG recommended choices which
tend to be closer to $0.1$ or smaller. We will see this trend continue in the examples that follow.
We have even observed higher tolerances (e.g., $\theta = .64$) performing well, though
larger $\theta$'s give rise to high AMG operator complexities. We note that
the operator complexities (not shown separately in the figure) range  from $1.11$ to $1.44$ for $\theta=.32$ with all three of these
 algorithm scenarios.
As smoothed aggregation AMG normally has very low complexities (e.g., near $1.1$ or lower) when compared
to classical AMG, it is often possible to more aggressively classify off-diagonals as weak (than might be
typical for classical AMG) and still end up with a practical method. 

Overall, the two rightmost columns (signed classical AMG scaling with traditional classification
versus symmetric SA scaling with cut-drop) are the best methods and perform
fairly similarly with cut-drop
being slightly better. In particular, there are a few interspersed lighter blue colors and one orange dot
with signed classical AMG scaling for larger $\theta$'s.
The interspersed colors and distributed lumping
are discussed more in Section~\ref{sec:results_3d}.
Finally, we note that
standard SA scaling and classification with the distance Laplacian (the middle column) does fairly well performing best
for $\theta = .08$.  Results for distributed lumping can be found in appendix Figure~\ref{fig:live2dDistribLump}.

\subsection{3D Tensor-Product Stretched
  Poisson}\label{sec:results_3d}
3D meshes 
are generated in a similar fashion to the Section~\ref{sec:results_2d} example
except for the addition of a fixed $z$ dimension containing 80 cells.
That is, twenty logarithmically spaced stretch factors
between 0.5 and 200 are investigated for the $x$ and $y$ dimensions.
Dirichlet conditions are
applied on the $y = 0$ boundary while  Neumann conditions are applied on other boundaries.
CG is again employed,
stopping when the residual is reduced by 10 orders of magnitude.
\begin{figure}[h!]\label{fig:live3c}
\hskip -.05in \scalebox{.215}[.215]{\includegraphics[trim = 3.4in 1.9in 2.0in 1.0in, clip            ]{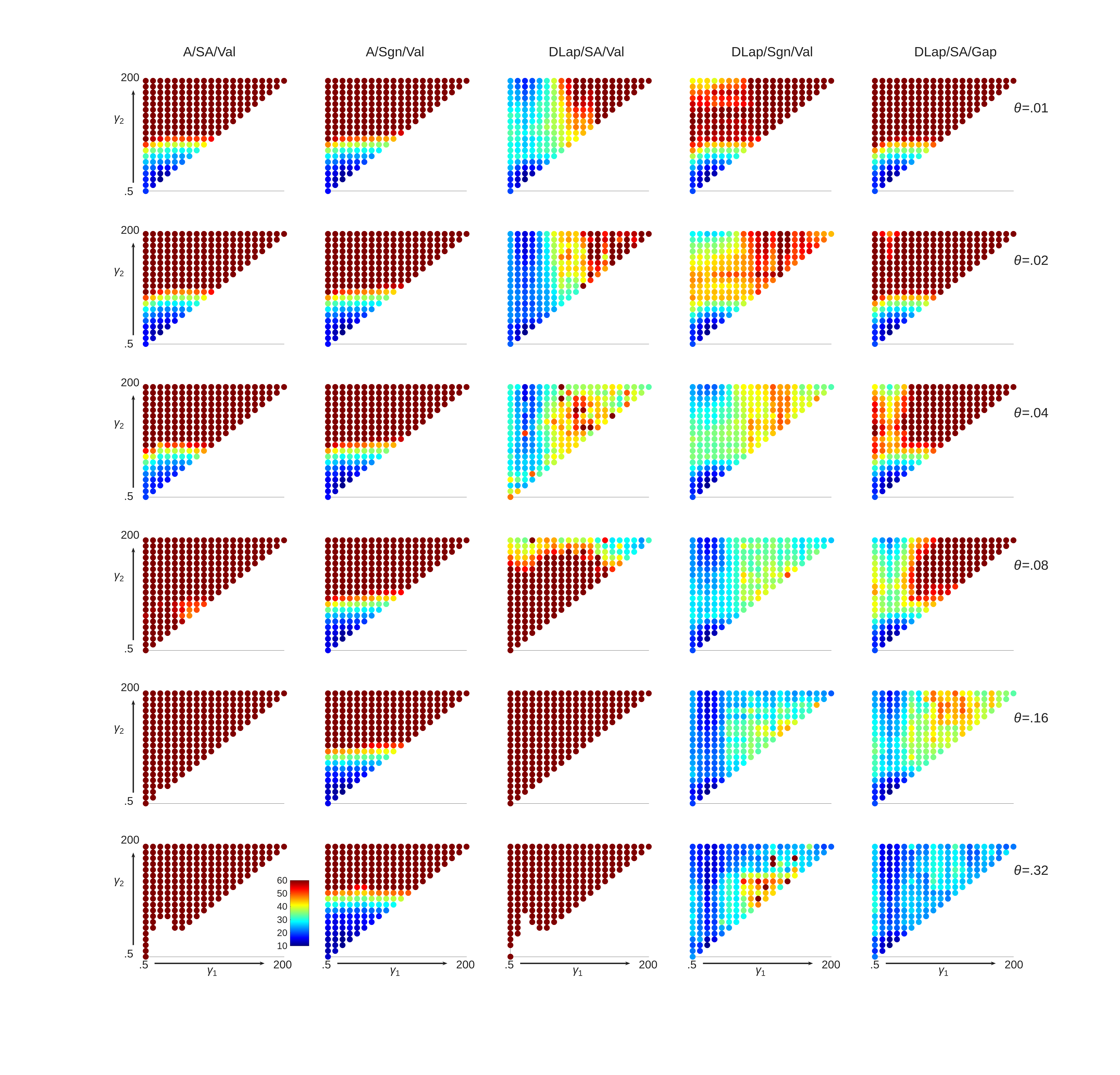}}
\vspace*{-.3in}
\caption {Iterations{\hskip .01in}$\times${\hskip .01in}operator{\hskip-.02in} complexity, with the color scale from 10 to 60 for the 3D mesh
Columns represent algorithm choices (described in
Section~\ref{sec:results_2d}) and rows represent tolerances,
$\theta$. Diagonal lumping is use for all experiments.}
\end{figure}
Figure~\ref{fig:live3c} shows traditional diagonal lumped results.
Generally, it can be seen that the distance Laplacian algorithms tend to
perform better over a range of $\theta$'s.  {\sf DLap/Sgn/Val} works well over the largest
range of $\theta$'s while the best individual case occurs for {\sf Dlap/SA/Gap} with $\theta = .32$.
\REMOVE{
(column 4) seems somewhat less sensitive to the $\theta$ choice, though there is no $\theta$ value
as good as the $\theta = .64$ {\sf A/Sgn/Val} or  $\theta = .32$ {\sf DLap/Sgn/Gap} (cut-drop) runs.
Clearly, the traditional SA scheme (column 1) is the worst performing. The traditional SA scaling and classification
scheme (column 3) performs somewhat better with the distance Laplacian, but it is the worst performing distance Laplacian,
even though it was reasonably competitive for the 2D example. Further, the best $\theta$ ($= .08$) for the 2D case is
not effective at all for the 3D case.
}

\begin{figure}[h!]\label{fig:live3c_lump}
\hskip -.05in \scalebox{.215}[.215]{\includegraphics[trim = 3.4in 1.9in 2.0in 1.0in, clip            ]{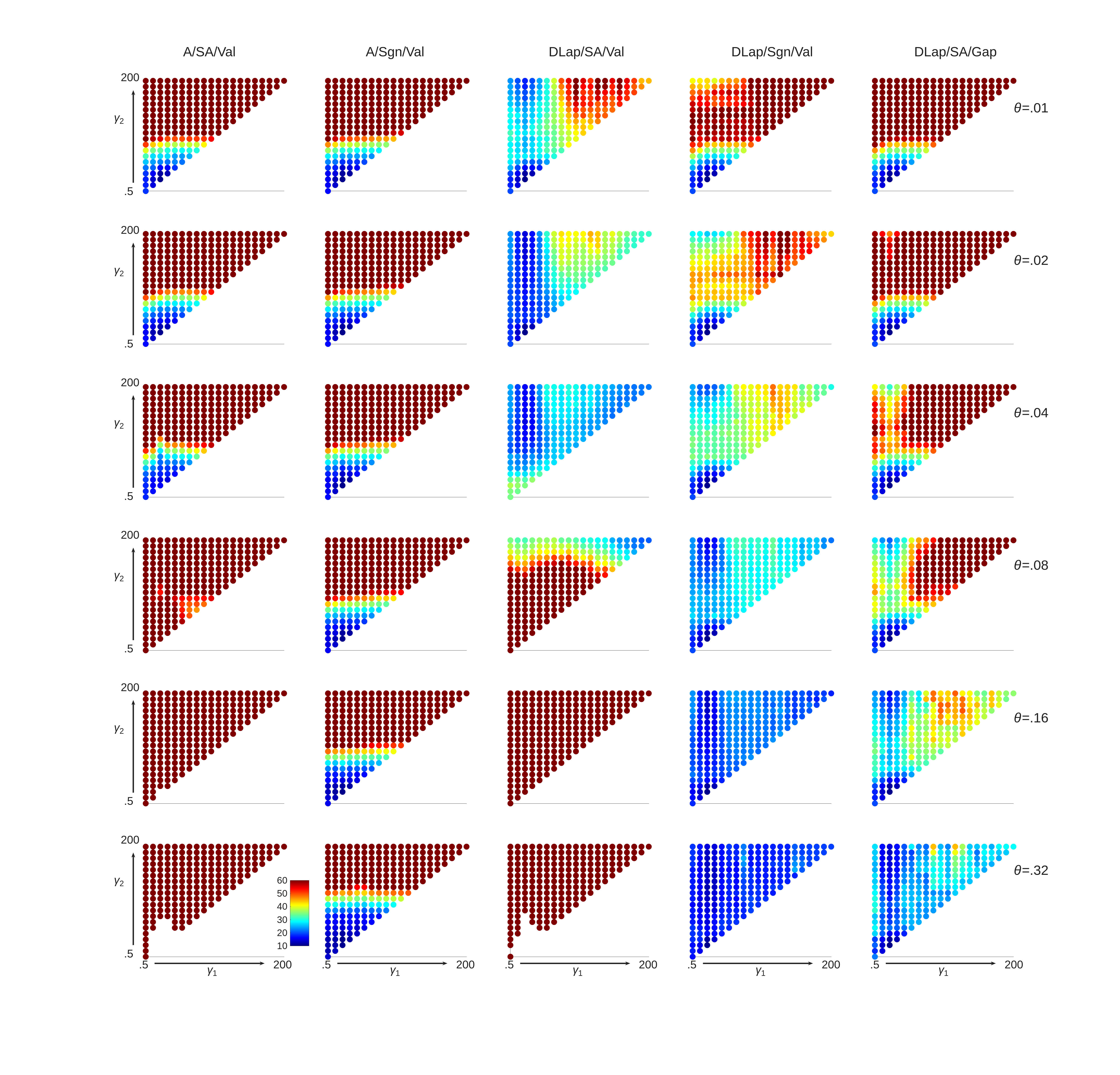}}
\vspace*{-.3in}
\caption{Iterations{\hskip .01in}$\times${\hskip .01in}operator{\hskip-.02in} complexity, with the color scale from 10 to 60 for the 2D mesh
Columns represent algorithm choices (described in
Section~\ref{sec:results_3d}) and rows represent tolerances,
$\theta$. Distributed lumping is use for all experiments.}
\end{figure}
Figure~\ref{fig:live3c_lump} shows the same scenarios using distributed lumping. Here we can see that
distributed lumping often modestly improves results. 
In general, the biggest impact occurs when a diagonal entry
in the dropped matrix became negative due to the traditional lumping procedure (which happens
more frequently at higher thresholds).
Overall, the biggest improvement is seen for {\sf DLap/Sgn/Val}
which is clearly the best performing combination.
The method now performs well for $\theta = .08, .16, ~\mbox{and}~ .32$.
These were also good $\theta$ choices for the 2D examples, where distributed lumping (see appendix Figure~\ref{fig:live2dDistribLump})
also removes the speckled dots seen in Figure~\ref{fig:live2d} for {\sf DLap/Sgn/Val}.
The  {\sf DLap/SA/Val} results are also noticeably
improved for $\theta = .02 ~\mbox{and}~ .08.$ Notice that the best $\theta$'s for {\sf DLap/SA/Val} are now somewhat smaller
than for the 2D problem. In general, we find guessing an appropriate $\theta$ choice for {\sf DLap/SA/Val} can be challenging
compared with choosing $\theta$ for {\sf DLap/Sgn/Val} or {\sf DLap/SA/Gap}.
Finally, we note that distributed lumping is not guaranteed to improve AMG convergence. We see that in many
cases, the distributed lumping and diagonal lumping give similar results. This can occur when standard lumping is not
problematic (i.e., it does not reduce diagonal entries by {\it too much}) or when convergence limitations are due
to sub-optimal classification choices so that changes to the lumping procedure are not impactful. This later case
occurs in many of the sub-plots with a large number of red dots.
As one can see, there are also situations
where standard lumping might perform better. In several of these cases, we've observed that standard lumping produces
a negative diagonal entry leading to a large negative eigenvalue estimate. This in turn results in a very small magnitude
(and negative) prolongator smoothing damping parameter $\omega$ in \eqref{eq: prolong smooth}, which oddly enough can have a fortunate beneficial effect. In particular,
the resulting grid transfer operator resembles piecewise-constant interpolation. 
While this interpolation is sub-optimal (when compared to smoothed grid transfer functions), it
may actually be better than that obtained with large $\omega$ prolongator smoothing
if the the coarse discretization operator is flawed due to earlier poor classification choices.
The {\sf DLap/SA/Gap} with $\theta = .32$ case is where some of the distributed lumping results seem a bit inferior.
We have noticed that with large thresholds on coarser grids (where some fill-in has occurred), {\sf DLap/SA/Gap} might classify many
small scaled SOC matrix entries as strong when there is not a nice gap in the magnitude of values within some rows.

\REMOVE{
It is interesting to note that {\sf A/Sgn/Val},
{\sf Dlap/SA/Gap} all do well with
with the distance Laplacian, which now appears to be the best overall choice in that it performs
fairly well over a range of $\theta$'s and for $\theta = .16$ or $\theta = .32$ leads to fairly robust performance.
For the remainder of the paper, we give only results with distributed lumping as it generally improves
the convergence behavior over traditional diagonal lumping.
A quantitative comparison of algorithm convergence
can be found in Table~\ref{tbl:live_3c-02262025}.
}

\REMOVE {

The traditional SA choices (leftmost column)

We note that dLapSa generally
does poorly on problems when the stretches approach 10:1
and that solver construction failures begin to appear
for $\theta \ge .08$ (seen as missing dots on the plot).
discClass once again does well for large $\theta$ but its performance
degrades for smaller thresholds on heavily stretched meshes.
By contrast,
dLapCut has consistently small iteration counts for any
choice of $\theta\in[.08,0.32]$.  Results become much
more erratic for dLapCut when $\theta=0.64$, but we see no cases of
solver failure.
}

\subsection{3D Radial Trisection Stretched Poisson}\label{sec:results_3d_radtri}
\begin{figure}[h!]\label{fig:radial}
\vskip -.1in
\begin{center}
 \includegraphics[trim = 5.8in 1.25in 5.65in 2.840in, clip,scale=0.06]{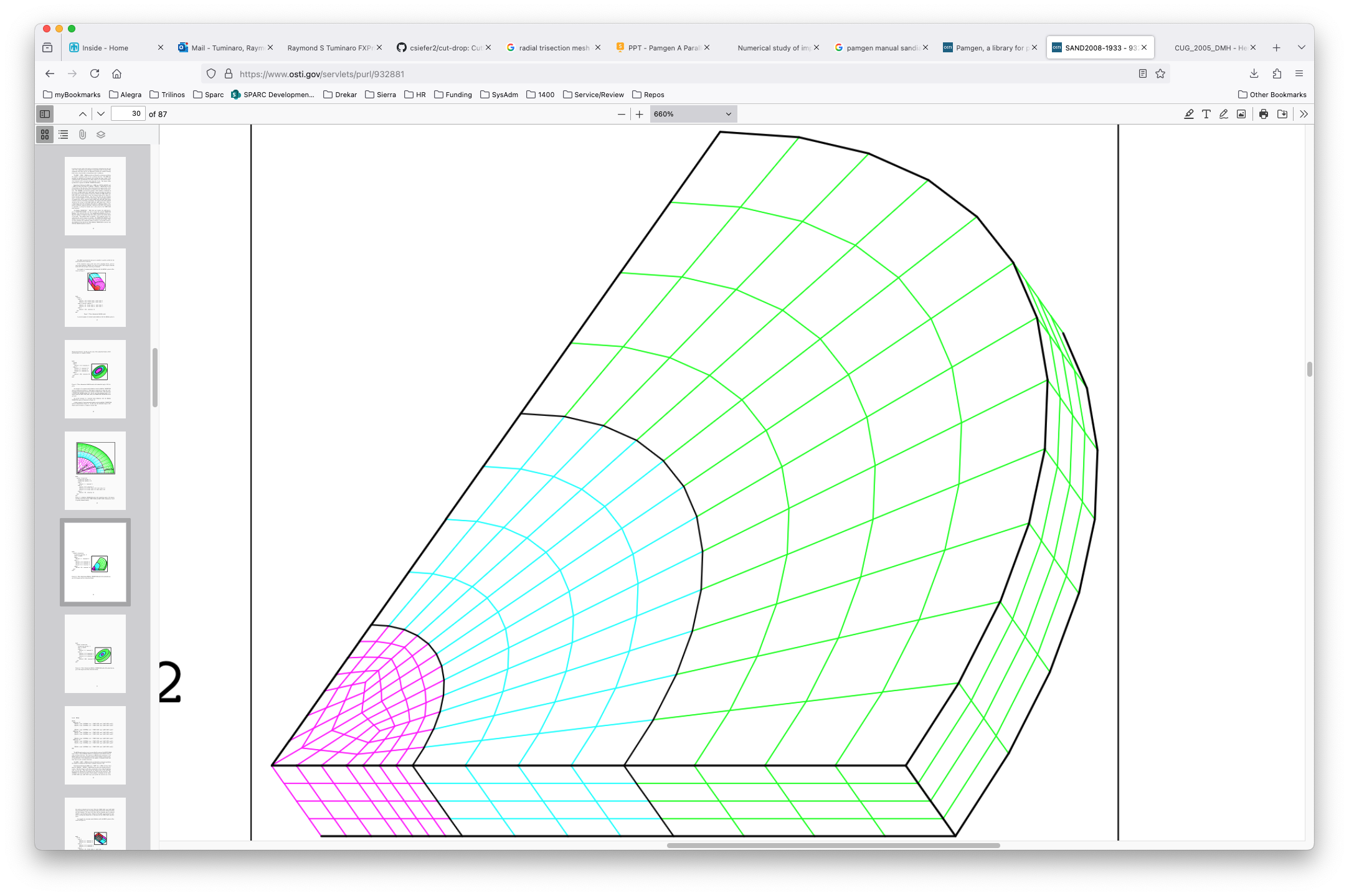}
\end{center}
\vspace*{-.1in}
\caption{Sample radial trisection mesh.}
\end{figure}

For some applications is not uncommon to pair a
radial mesh with a less structured mesh near the origin
to avoid degenerate elements, as in (e.g., Figure~\ref{fig:radial}).
We vary these innermost-region transition blocks,
which we call \textit{trisection blocks}, with values 1, 2, 3, 4 and 6.
Depending on the number of trisections, we start by
stretching the radial region by a factor of 1.0658 to 1.6426 and
increase stretching logarithmically up to 100 leading
to 205 cases.
Dirichlet conditions 
are applied on the outer radius while other boundaries have Neumann conditions.

\begin{figure}[h!]\label{fig:live4b}
\vskip -.1in
\begin{center}
\hskip -.05in
\scalebox{.175}[.125]{\includegraphics[trim = 3.4in 1.9in 0.0in 0.9in, clip ]{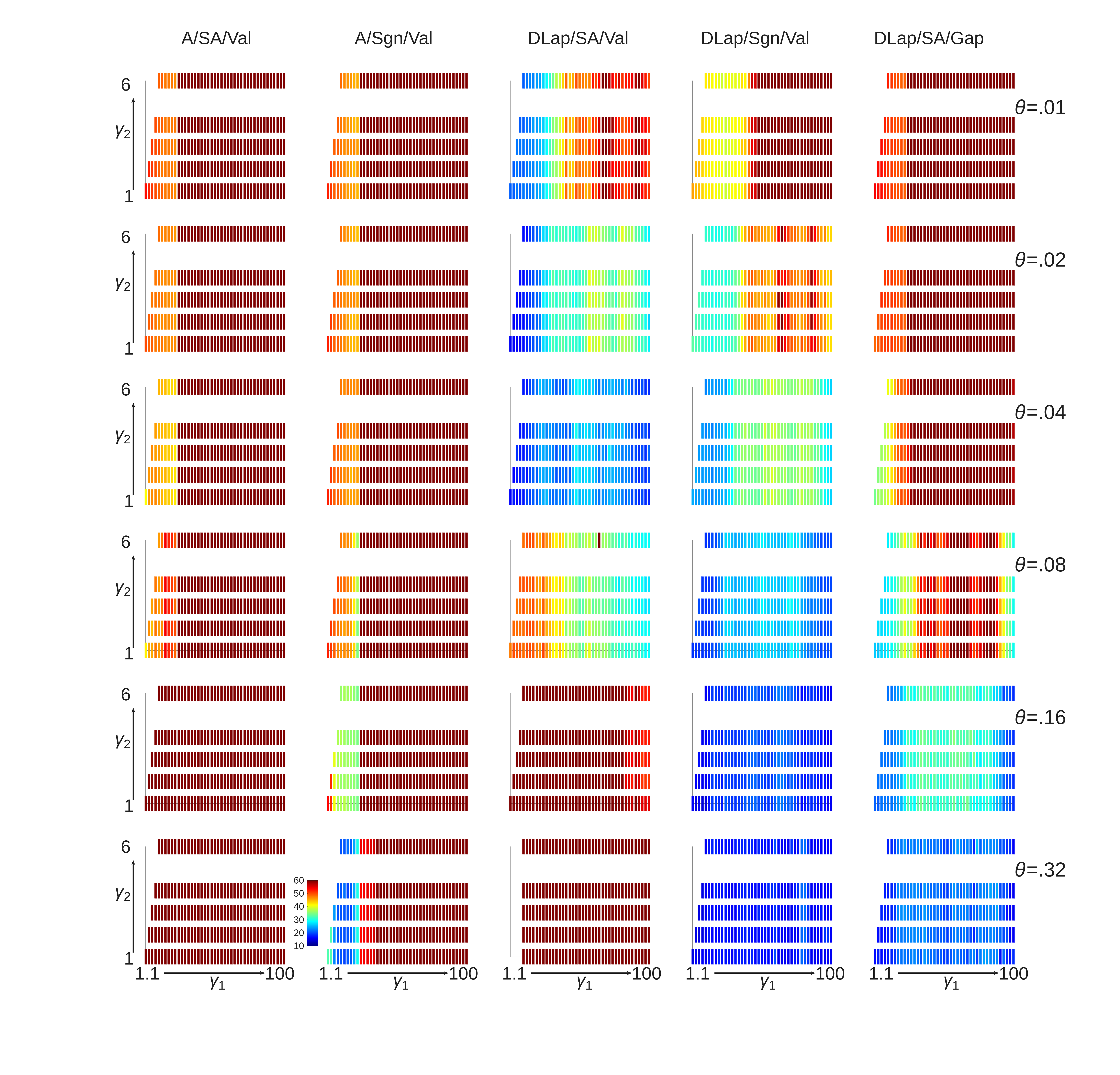}}
\end{center}
\vspace*{-.2in}
\caption{Iterations~$\times$~operator{\hskip -.02in} complexity (distributed lumping) ~for radial hexahedron case.}
\end{figure}
Figure~\ref{fig:live4b} shows CG results with distributed lumping where now each run is plotted as a vertical bar.
Convergence is declared when the residual is reduced by 10 orders of magnitude.
The number of trisection blocks $\gamma_2$ takes on the values of $1, 2, 3, 4,$ and $6$
(the gap between the top and second line in the sub-figures is due to the omission of $\gamma_2 = 5$).
These results are fairly consistent with those from the 3D tensor mesh example.
We again see that the distance Laplacian methods outperform those that use $A$ for the SOC matrix.
Traditional SA ({\sf A/SA/Val} in column 1) is the worst. Classical AMG scaling (column 2) improves the
situation, but relies heavily on automatically discarding all positive off-diagonals.
Among distance Laplacian versions {\sf DLap/Sgn/Val} (column 4) is the best followed by
{\sf DLap/SA/GAP} (column 5) and then {\sf DLap/SA/Val} (column 3).

\subsection{2D Boundary Layer Meshes with Circular Cavities }\label{sec:results_2d_cavities}
\begin{figure}\label{fig:circular cavities}
\includegraphics[trim = 0.05in 0.15in 0.75in 0.65in, clip, width= 6.0in,height= 2.5in ]{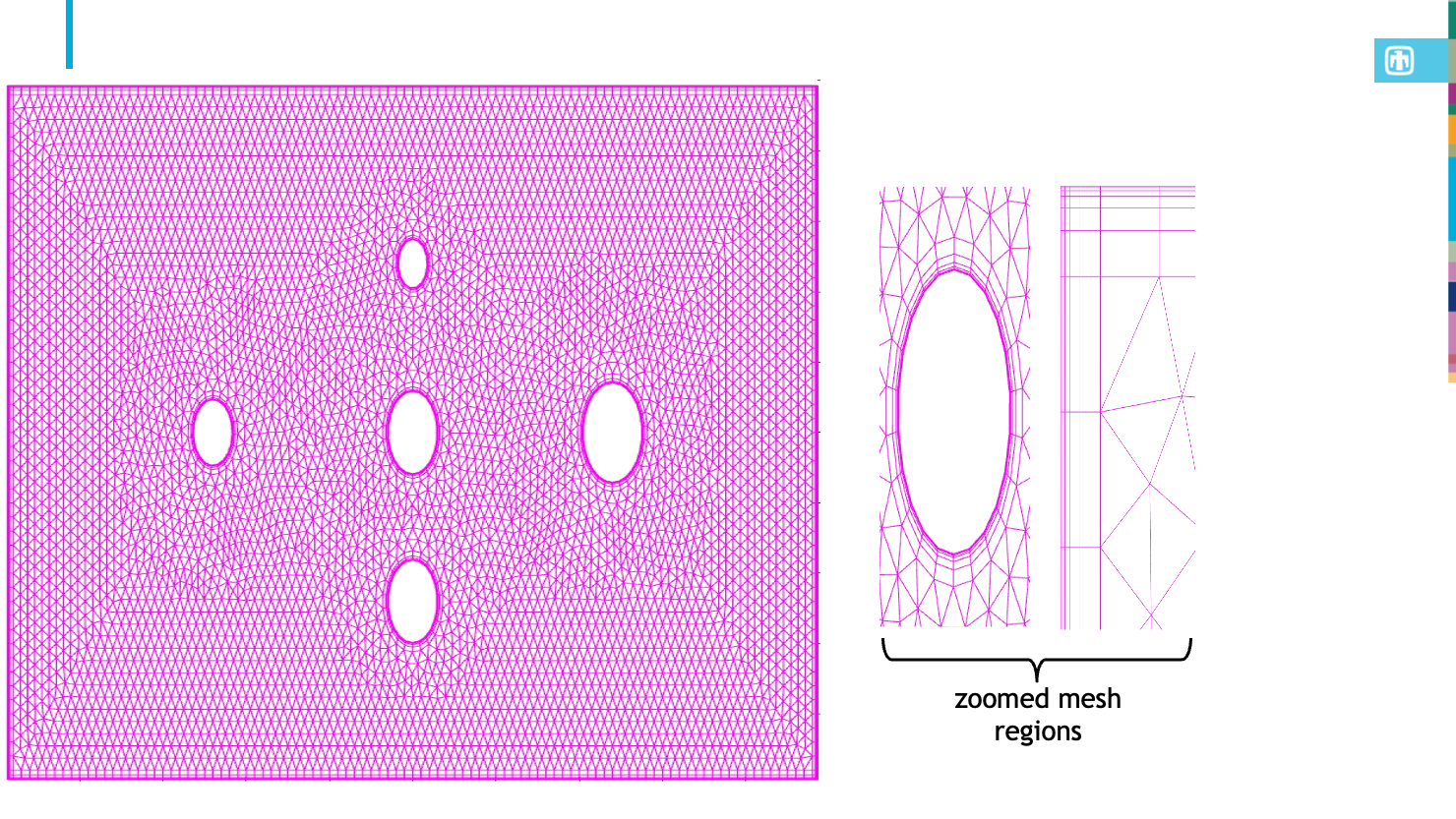}
\caption{Five cavity mesh.}
\end{figure}

We now consider a mesh with  circular cavities and includes boundary layers (Figure~\ref{fig:circular cavities}).
Dirichlet boundary conditions are present on the top and bottom
domain boundaries as well as on the central cavity boundary.
The top Dirichlet condition enforces the solution to be one while the other Dirichlet boundaries enforce
a zero solution.  All other boundary conditions are homogeneous Neumann.
Triangular elements are used for non-boundary layer portions of the mesh
while quadrilateral elements are used for the boundary layers.  The CG algorithm is used
with a relative residual tolerance of $10^{-10}$. The only nonzero right hand side entries
correspond to the top Dirichlet boundary.
\begin{table}[h!]\begin{center}\begingroup \setlength{\tabcolsep}{4pt}\footnotesize\begin{tabular}{|r||ccc|ccc|ccc|ccc|ccc|}
\hline
$\theta$&\multicolumn{3}{c|}{A/SA/Val}&\multicolumn{3}{c|}{A/Sng/Val}&\multicolumn{3}{c|}{DLap/SA/Val}&\multicolumn{3}{c|}{DLap/Sng/Val}&\multicolumn{3}{c|}{DLap/SA/Gap} \\
        & SL & $\mathcal{D}$L &$\alpha$& SL  & $\mathcal{D}$L & $\alpha$ & SL & $\mathcal{D}$L & $\alpha$ & SL & $\mathcal{D}$L & $\alpha$ & SL & $\mathcal{D}$L & $\alpha$ \\
\hline
    .01   &89 & 88 &1.06  & 73 & 73 & 1.07  & 22 &17 &1.07  &  27 & 25 & 1.08   & 32 & 32 & 1.08 \\
    .01   &90 & 90 &1.11  & 80 & 80 & 1.29  & 36 &36 &1.17  &  50 & 49 & 1.17   & 54 & 54 & 1.17 \\
    .02   & 92& 92 &1.11  & 80 & 80 & 1.29  & 33 &33 &1.16  &  45 & 45 & 1.16   & 50 & 50 & 1.18 \\
    .04   & 88& 88 &1.11  & 80 & 80 & 1.29  & 23 &23 &1.20  &  41 & 41 & 1.18   & 45 & 45 & 1.16 \\

    .08   & 81& 81 &1.16  & 79 & 79 & 1.29  & 31 &31 &1.36  &  40 & 40 & 1.18   & 45 & 45 & 1.17 \\
    .16   & 94& 93 &1.35  & 72 & 72 & 1.31  & 75 &75 &1.40  &  35 & 35 & 1.18   & 44 & 44 & 1.17 \\
    .32   & - & -  &      & 34 & 34 & 1.37  & - & -  &      &  24 & 24 & 1.28   & 32 & 32 & 1.32 \\
\hline
\end{tabular}
\endgroup
\caption{Iteration counts for standard lumping SL and distributed lumping $\mathcal{D}$L  and AMG operator complexities $\alpha$ for five cavity mesh.
} \label{tbl:cavities}
\end{center}
\end{table}
Table~\ref{tbl:cavities} gives the results.
\REMOVE{
\begin{table}\begin{center}
\begin{tabular}{|c||c|c|}
\hline
Threshold & Classic & Cut-Drop \\
\hline\hline
.9       &    f    & 10/1.5   \\
.75      &    f    & 14/1.3   \\
.5       &    f    & 14/1.3   \\
.3536    &    f    & 19/1.3   \\
.25      & 60/1.2  & 20/1.2   \\
.1768    & 49/1.3  & 23/1.2   \\
.125     & 17/1.2  & 21/1.2   \\
.0884    & 16/1.2  & 28/1.2   \\
\hline
\end{tabular} ~~~~~~~~~~~~~~
\begin{tabular}{|c||c|c|}
\hline
Threshold & Classic & Cut-Drop \\
\hline\hline
.0625    & 20/1.2  & 24/1.2   \\
.0442    & 21/1.2  & 25/1.2   \\
.0312    & 20/1.2  & 28/1.2   \\
.0221    & 22/1.1  & 23/1.2   \\
.0156    & 26/1.1  & 24/1.2   \\
.0110    & 43/1.1  & 25/1.2   \\
.0078    & 49/1.1  & 28/1.2   \\
\hline
\end{tabular}
\caption{Iteration counts/complexities  for 2D boundary layer cavity mesh.'f' indicates a failure to solve
the problem either because of a breakdown in the setup phase (e.g., all entries dropped from a row) or not satisfying the convergence criteria in 100 iterations. }\label{tbl:cavities}
\end{center}
\end{table}
}
A dash entry indicates that the solver did not converge within 200 iterations.
These results mirror those of the previous problems. In particular, the distance Laplacian variants generally perform better with
{\sf DLap/Sng/Val} being the best while {\sf DLap/SA/Gap} is the second best. For the cavity problem, there is no significant
difference between using standard and distributed lumping.
\REMOVE {
using
While both algorithm are able to solve the problem effectively (e.g., 20 iterations), the range
of successful thresholds for the cut-drop algorithm is much greater than for the classic algorithm.
Convergence is achieved even with a threshold of $.9$ requiring only 10 iterations, though the
AMG operator complexity is a bit higher than the other cases.
It is worth noting that with a threshold of $0.0$, $49$ iterations are required and the complexity is $1.1$.
Additionally, we remark that using $A$ for the SOC matrix does quite poorly often leading to
little improvement over the non-dropping algorithm.
}

\subsection{Meshes from Applications}\label{sec:application meshes}
\begin{figure}[h!]\label{fig:plumeMesh}
\begin{center}
\includegraphics[trim = 1.5in 0.0in 0.0in 0.0in, clip,scale=0.171]{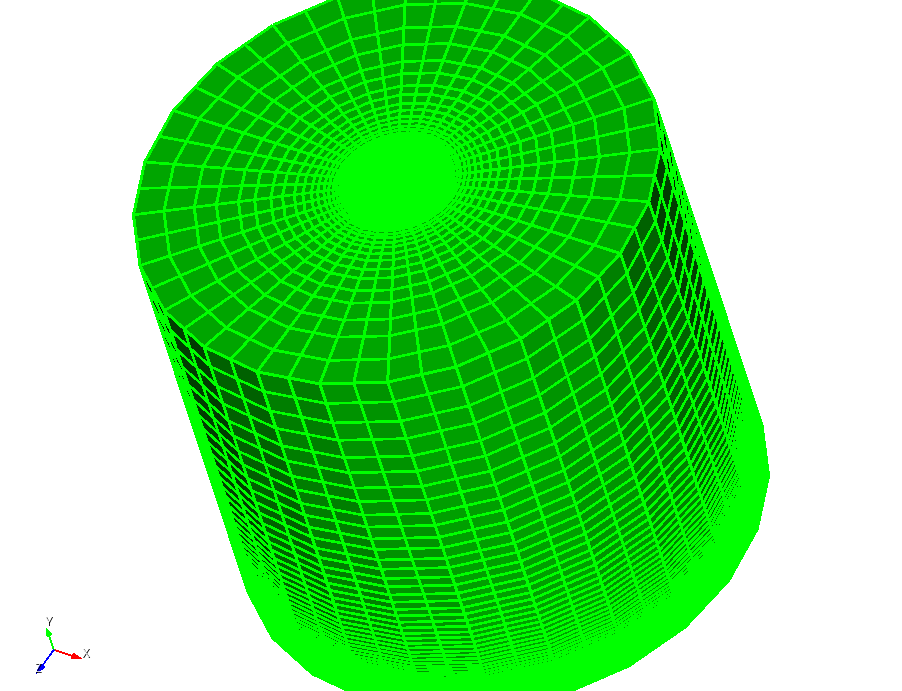}
~~~
\includegraphics[trim = 0.0in 0.0in 0.0in 0.0in, clip,scale=0.3344]{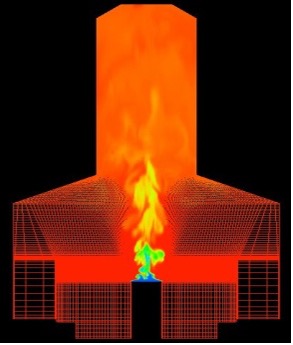}
\end{center}
\vskip -.1in
\caption{Partially zoomed image of simplified mesh used for Table~\ref{tbl:plumeTable} (left) and mesh
         densities with gas density contours overlaid for actual helium plume studies (right).}
\end{figure}

We begin this sub-section with a cylinder mesh (left side of Figure~\ref{fig:plumeMesh}) that has been used as a simple benchmark
for more sophisticated helium plume meshes and studies~\cite{plume}.
All results in this sub-section investigate Poisson problems on stretched meshes that are
discretized using linear first order hexahedral finite elements.  In all of these examples,
Neumann boundary conditions are applied on all but one surface which instead has Dirichlet boundary conditions.
For the cylinder mesh, the Dirichlet boundary condition coincides with one of the circular surfaces
and the matrix $A$ has size $772,585 \times 772,585 $.
To limit the required computations needed for the generally larger and more complex meshes in this sub-section,
we now employ AMG preconditioned GMRES(300) and convergence is declared when the two-norm of the initial residual is reduced by $10^{-6}$.
This change was primarily motivated by our last example, which needed to be run in serial due to some mesh partitioning and format conversion limitations.

Table~\ref{tbl:plumeTable} gives iteration counts and AMG operator complexities.
\begin{table}[h!]\begin{center}\begingroup \setlength{\tabcolsep}{4pt}\footnotesize\begin{tabular}{|r||ccc|ccc|ccc|ccc|ccc|}
\hline
$\theta$&\multicolumn{3}{c|}{A/SA/Val}&\multicolumn{3}{c|}{A/Sng/Val}&\multicolumn{3}{c|}{DLap/SA/Val}&\multicolumn{3}{c|}{DLap/Sng/Val}&\multicolumn{3}{c|}{DLap/SA/Gap} \\
        & SL & $\mathcal{D}$L &$\alpha$& SL  & $\mathcal{D}$L & $\alpha$ & SL & $\mathcal{D}$L & $\alpha$ & SL & $\mathcal{D}$L & $\alpha$ & SL & $\mathcal{D}$L & $\alpha$ \\
\hline
    .01   & 52& 52 &1.06  & 44 & 44 & 1.07  & 19 & 15 &1.07  &  21 & 21 & 1.08   & 27 & 27 & 1.08 \\
    .02   & 54& 54 &1.05  & 43 & 43 & 1.07  & 17 & 12 &1.11  &  17 & 17 & 1.06   & 23 & 23 & 1.09 \\
    .04   & 46& 45 &1.08  & 43 & 43 & 1.07  & 19 & 15 &1.22  &  19 & 14 & 1.08   & 23 & 23 & 1.06 \\
    .08   & 50& 45 &1.17  & 33 & 33 & 1.10  & 25 & 23 &1.32  &  17 & 13 & 1.11   & 19 & 19 & 1.08 \\
    .16   & 80& 50 &1.24  & 33 & 33 & 1.13  & 49 & 49 &1.19  &  19 & 10 & 1.18   & 15 & 15 & 1.13 \\
    .32   & 65& 64 &1.11  & 20 & 20 & 1.24  & 85 & 85 &1.06 &   18 &  7 & 1.27   & 12 & 11 & 1.20 \\
\hline
\end{tabular}
\endgroup
\caption{Iteration counts and AMG operator complexities $\alpha$ for standard lumping SL and distributed lumping $\mathcal{D}$L
on cylinder mesh problem.  } \label{tbl:plumeTable}
\end{center}
\end{table}
The AMG operator complexities for standard and distributed lumping differ by at most $.01$. To reduce
clutter, only the largest of these two complexities is shown for each dropping algorithm/$\theta$ choice in all tables within this sub-section.
For all runs, the number of AMG levels varies between $3$ and $5$.

From the results, it is clear that variants employing the distance Laplacian SOC matrix are generally better
than those that employ $A$ as SOC matrix. The overall best method is {\sf DLap/Sgn/Val} in that it successfully
solves the problem in under $15$ iterations with distributed lumping over a wide range of $\theta$ values. We
can also see that distributed lumping generally either reduces the iteration count or has the same iteration
count as standard diagonal lumping.
In one case, distributed
lumping reduced the iteration count from $ 80$ to $ 50$ and
distribution lumping requires about half as many iterations as standard lumping for the best performing combination:
{\sf DLap/Sgn/Val} with $\theta \ge .16$.
Finally, we note that for realistic helium plume simulations depicted on the right side of Figure~\ref{fig:plumeMesh},
a sophisticated control volume discretization technique is instead employed
and the solution evolution requires repeated Poisson solves in a pressure-projection type of method.

\begin{figure}[h!]\label{fig:jetInCrossflow}
\begin{center}
\includegraphics[trim = 1.5in 0.0in 0.0in 0.5in, clip,scale=0.171]{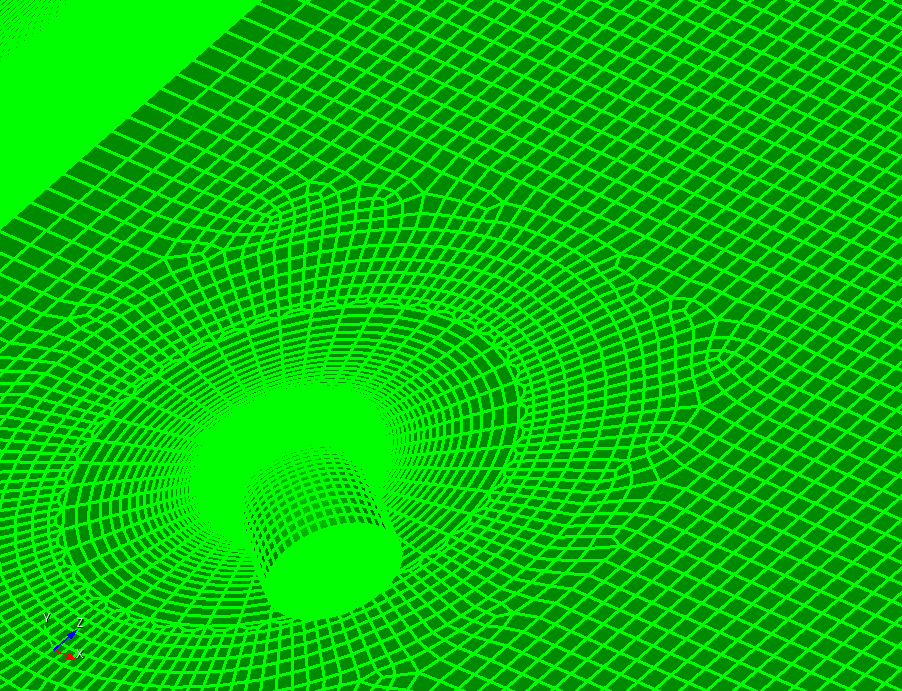}
\end{center}
\vskip -.1in
\caption{Zoomed image of mesh used for Table~\ref{tbl:jetInCrossflow} focusing on inlet where jet enters cross flow.}
\end{figure}

We next consider a mesh that is normally used to simulate a jet of fluid perpendicularly entering a cross flow in \textsc{Nalu}~\cite{nalu},  requiring refined meshes 
where the jet and
cross flow meet to resolve the complex interactions. Figure~\ref{fig:jetInCrossflow} illustrates a highly zoomed portion of the mesh which essentially
corresponds to a small inlet region that is connected to a much larger box near one of its corners.
The discrete matrix $A$ is $3,548,688 \times 3,548,688$ and the one Dirichlet surface coincides with
the circular face at the tip of the inlet.
Table~\ref{tbl:jetInCrossflow} gives the iteration counts and AMG operator complexities.
\begin{table}[h!]\begin{center}\begingroup \setlength{\tabcolsep}{4pt}\footnotesize\begin{tabular}{|r||ccc|ccc|ccc|ccc|ccc|}
\hline
$\theta$&\multicolumn{3}{c|}{A/SA/Val}&\multicolumn{3}{c|}{A/Sng/Val}&\multicolumn{3}{c|}{DLap/SA/Val}&\multicolumn{3}{c|}{DLap/Sng/Val}&\multicolumn{3}{c|}{DLap/SA/Gap} \\
        & SL & $\mathcal{D}$L &$\alpha$& SL  & $\mathcal{D}$L & $\alpha$ & SL & $\mathcal{D}$L & $\alpha$ & SL & $\mathcal{D}$L & $\alpha$ & SL & $\mathcal{D}$L & $\alpha$ \\
\hline
    .01  &  51  &  51 &  1.04   &   41  &   41 & 1.05    &   26  &  22 & 1.21  &  51  &  29 & 1.16  &  42  &  42 & 1.15 \\
    .02  &  50  &  50 &  1.04   &   41  &   41 & 1.06    &   81  &  32 & 1.27  &  39  &  25 & 1.19  &  36  &  36 & 1.18 \\
    .04  &  48  &  48 &  1.10   &   40  &   40 & 1.06    &  103  &  43 & 1.39  &  45  &  26 & 1.23  &  31  &  31 & 1.22 \\
    .08  & 300  & 275 &  1.19   &   35  &   35 & 1.12    &  123  &  88 & 1.48  &  43  &  19 & 1.29  &  31  &  31 & 1.26 \\
    .16  & 497  & 377 &  1.32   &   33  &   33 & 1.16    &  166  & 160 & 1.36  &  44  &  12 & 1.35  &  26  &  26 & 1.32 \\
    .32  & 381  & 357 &  1.28   &   17  &   18 & 1.35    &  263  & 263 & 1.22  &  50  &   9 & 1.42  &  15  &  15 &  1.37 \\
\hline
\end{tabular}
\endgroup
\caption{Iteration counts and AMG operator complexities $\alpha$ for standard lumping SL and distributed lumping $\mathcal{D}$L on
inlet jet problem.
} \label{tbl:jetInCrossflow}
\end{center}
\end{table}
We see that using the distance Laplacian as the SOC matrix generally improves the convergence rate when compared to its counterpart
using the matrix $A$ as the SOC matrix. However, it should be noted that {\sf A/Sgn/Val} is better than {\sf DLap/Sgn/Val} when
using standard diagonal lumping, but the converse is true when using distributed lumping (which is the best overall performing method).

The results for the jet inlet problem mirror those of the cylinder mesh in that the overall best method is {\sf DLap/Sgn/Val} with distributed lumping,
which provides converged solutions for this problem in under
$20$ iterations with $\theta \ge .08$.
Additionally, we see that the iteration counts are not too sensitive for the algorithm combinations that use either signed classical scaling or cut-drop
({\sf A/Sng/Val}, {\sf DLap/Sng/Val} or {\sf DLap/SA/Gap}) when compared to the other two algorithm combinations
({\sf A/SA/Val} or {\sf DLap/SA/Val})
While distributed lumping {\sf DLap/SA/Val} runs are somewhat competitive, it is clearly the worst distance Laplacian scheme.
Similar to the 3D stretched tensor meshes, {\sf DLap/SA/Gap} generally performs better than {\sf DLap/Sgn/Val} with standard
lumping, but the opposite is true once distributed lumping is introduced. That is, distributed lumping appears to have a more positive impact on
the distance Laplacian signed classical scaling method than on the cut-drop method.

\begin{figure}
~~~~~~~~~~\includegraphics[height = 3.8cm,width = 4.8cm]{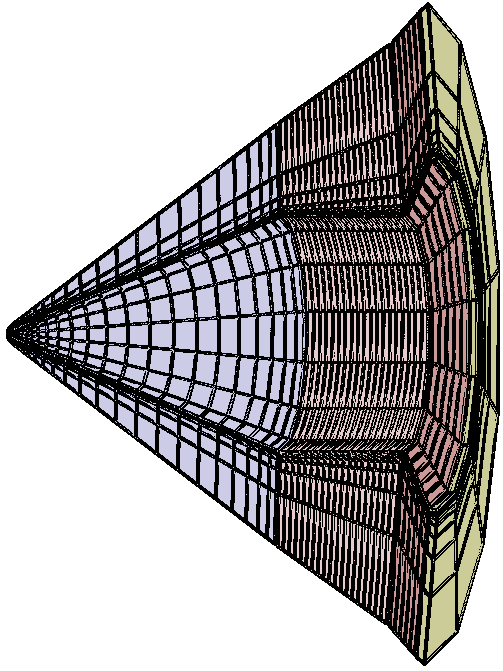}~~~~~~~~~~~~~~~
\includegraphics[trim = 1.5in 0.0in 2.0in 0.0in, clip, height = 3.8cm, width = 3.8cm]{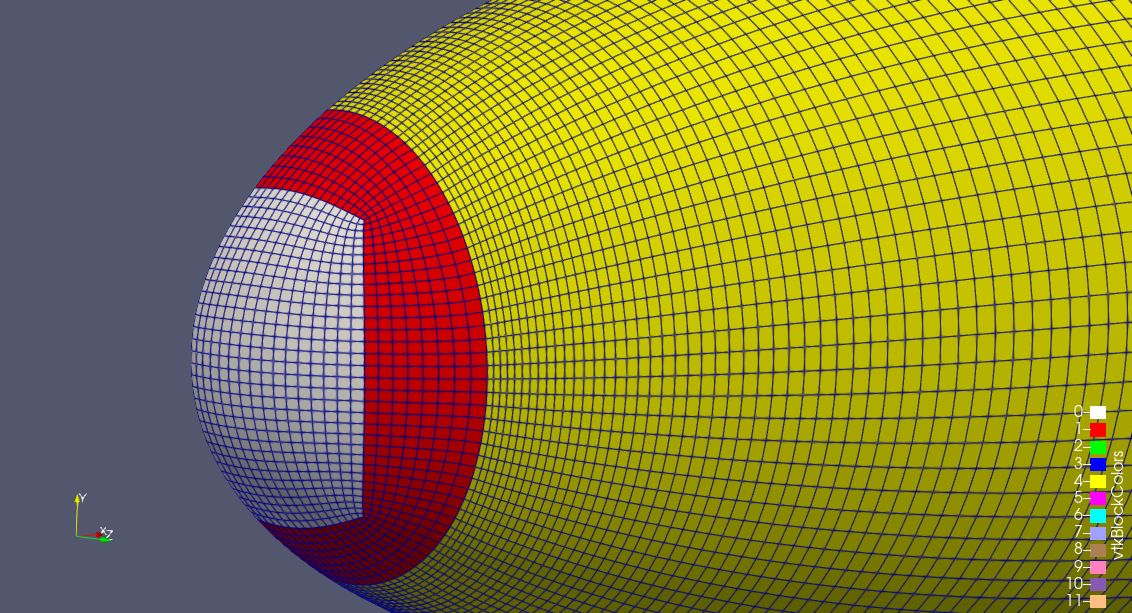}
\hfill~
\caption{
\emph{Left:} 6 block coarse version of HIFiRE-1 mesh. Each color denotes a different block structured mesh region.
\emph{Right:} Zoom of outer HIFiRE-1 mesh surface near nose.}
\label{fig:sparcMesh}
\end{figure}

We conclude this subsection with a mesh (see Figure~\ref{fig:sparcMesh}) originally generated with Pointwise~\cite{woeber2017mesh} that is
normally used for hypersonic flow of re-entry vehicles
in the Sandia Parallel Aerodynamics Reentry Code (\textsc{SPARC})~\cite{sparc}.
The computational domain corresponds to the space surrounding the vehicle with finer mesh resolution as one approaches
the vehicle surface (i.e., cylindrical cavity) and the nose cone to resolve boundary layers.
The discrete matrix $A$ is $4,298,409 \times 4,298,409$ and the Dirichlet surface corresponds to the
yellow horseshoe shaped face furthest from the cone.
Table~\ref{tbl:hifire} gives the iteration counts and AMG operator complexities.
\begin{table}[h!]\begin{center}\begingroup \setlength{\tabcolsep}{4pt}\footnotesize\begin{tabular}{|r||ccc|ccc|ccc|ccc|ccc|}
\hline
$\theta$&\multicolumn{3}{c|}{A/SA/Val}&\multicolumn{3}{c|}{A/Sng/Val}&\multicolumn{3}{c|}{DLap/SA/Val}&\multicolumn{3}{c|}{DLap/Sng/Val}&\multicolumn{3}{c|}{DLap/SA/Gap} \\
        & SL & $\mathcal{D}$L &$\alpha$& SL  & $\mathcal{D}$L & $\alpha$ & SL & $\mathcal{D}$L & $\alpha$ & SL & $\mathcal{D}$L & $\alpha$ & SL & $\mathcal{D}$L & $\alpha$ \\
\hline
    .01  & 413 & 413 & 1.04 & 348   & 348  & 1.07 &  92   &  43 & 1.42 &   34  & 33 & 2.28  & 101 & 101 & 2.16  \\
    .02  & 404 & 405 & 1.05 & 344   & 344  & 1.07 & 107   &  54 & 1.42 &   27  & 28 & 2.36  &  96 &  95 & 2.56  \\
    .04  & 373 & 370 & 1.08 & 344   & 344  & 1.07 & 246   &  65 & 1.45 &   26  & 22 & 2.14  &  86 &  84 & 2.54  \\
    .08  & 415 & 375 & 1.14 & 274   & 274  & 1.17 & 324   & 110 & 1.49 &   45  & 15 & 2.31  &  72 &  72 & 2.46  \\
    .16  &  -  &  -  &      & 271   & 271  & 1.20 & 215   & 191 & 1.45 &   57  & 11 & 2.37  &  57 &  57 & 2.40  \\
    .32  & 262 & 246 & 1.40 &  26   &  26  & 1.56 & 290   & 290 & 1.38 &   37  & 10 & 2.29  &  40 &  38 & 3.33  \\
\hline
\end{tabular}
\endgroup
\caption{Iteration counts and AMG operator complexities $\alpha$ for standard lumping SL and distributed lumping $\mathcal{D}$L
on HIFiRE-1 mesh problem.
} \label{tbl:hifire}
\end{center}
\end{table}
A dash indicates that a particular run did not converge within 500 iterations.
As with previous results, using $A$ as the SOC matrix
generally leads to poor convergence with the noticeable exception of {\sf A/Sgn/Val} which converges
well for a quite large $\theta = .32$ while converging poorly for other $\theta$'s.
The {\sf DLap/SA/Val} results with diagonal lumping are somewhat disappointing,
but they improve dramatically with distributed lumping. {\sf A/SA/Val}, {\sf A/Sgn/Val}, and {\sf DLap/SA/Val}
give by far the best operator complexities. Thus, if storage is a concern, {\sf A/Sgn/Val} with a
very high tolerance or {\sf DLap/SA/Val} with a low tolerance and distributed lumping
could be attractive.  We note that $\theta = .01$ is indeed a fairly good choice for {\sf DLap/SA/Val} and
distributed lumping as iterations counts rise with smaller $\theta = .001$.  (not shown here).
In terms of convergence over a range of $\theta$'s, however, the best method is
{\sf DLap/Sgn/Val} in conjunction with distributed lumping. In particular, with $\theta \ge .08$
the number of required iterations is nearly three times less than the best option with  {\sf DLap/SA/Val}
and is $2.5$ times better than {\sf A/Sgn/Val}. Although the its AMG operator complexity is $1.5 $ times higher than that
of  {\sf DLap/SA/Val}, it is robust over a range of $\theta$'s and does not rely on the sign
of the off-diagonal entries.
Finally,
we note that distributed lumping has a major impact for {\sf DLap/SA/Val} and {\sf DLap/Sgn/Val}, while having little influence on
the other three algorithm combinations.

\section{Conclusions and Future Work}\label{sec:conclusions}

We have studied strong and weak AMG classification when solving Poisson
problems on stretched meshes using linear finite elements.
We have centered our study on the choice of SOC matrix, the scaling
applied to the SOC matrix entries, and whether or not it is best
to classify based solely on a threshold or whether it can be advantageous
to instead look for gaps in the scaled SOC entries.
Additionally, we have defined a new distributed lumping procedure to account for
dropped/weak entries. Under modest assumptions, this new lumping procedure
is guaranteed to not reverse the sign of diagonal entries in the discretization matrix, which
could have serious detrimental convergence effects. Overall, we find
that the best algorithm combination for Poisson problems uses a distance Laplacian
operator to define the SOC matrix together with a non-symmetric
scaling scheme, which is normally used in classical AMG as opposed to smoothed aggregation
AMG. In conjunction with thresholds $\theta \ge .08$ and distributed lumping, this
combination performed consistently well on finite element discretized
Poisson problems defined over a wide range of meshes. Generally,
we found that fairly high thresholds also performed well for
{\sf DLap/SA/Gap}. While the best $\theta$ choice
for a more traditional smoothed aggregation criteria using either $A$
or a distance Laplacian SOC matrix could vary significantly.  The relatively consistent
best-$\theta$ choice for {\sf DLap/Sgn/Val} and {\sf DLap/SA/Gap}
is due to the way in which that they leverage the largest magnitude
off-diagonal entry of the distance Laplacian within each row. In particular, both
schemes guarantee that this entry is labeled strong along with other entries
whose magnitude is {\it close} to this largest-magnitude entry.
Overall, proper classification has a significant impact on AMG convergence
with the best performing variants requiring $10$'s of iterations while
poor performing variants might require $100$'s of iterations or not
even converge.

The main limitation in this study is that the distance Laplacian is mostly
appropriate for diffusion dominated systems without significant material
variations. An archive version of the paper {\it Smoothed aggregation algebraic multigrid for problems with heterogeneous and anisotropic materials},
by Firmbach, Phillips, Glusa, Popp, Siefert, and Mayr will appear shortly. This paper generalizes
the distance Laplacian idea to problems with material variation. A second paper is in the planning
phase which focuses on defining a SOC matrix by developing a crude approximation to $M^{-1} A$ where
$M$ is the finite element mass matrix and $A$ is the finite element stiffness matrix. These
approaches are somewhat complementary and are both intended to address limitations of the
distance Laplacian. We also plan to further adapt and evaluate the distributed
lumping scheme to more complex partial differential equations.

\REMOVE {

for a number of
linear finite
We systematically evaluated several

We introduced a new, cut-based variant of the smoothed aggregation distance
Laplacian dropping criterion, which does not use a single, global
threshold.  For both the cut-based algorithm and the traditional distance
Laplacian, we have calculated key parameter values which avoid
departure from brick coarsening for isotropic portions of the mesh.
Respecting this key value, we demonstrate that the cut-based criterion
consistently produces iteration counts comparable to or better than
the traditional dropping criterion on a variety of different problems
in both 2D and 3D.
Perhaps more importantly, the results from
the cut-based algorithm appear robust over a wide range of parameter
choices for the meshes tested.
While this study focused on the Poisson equation, we believe that this
technique could be directly applied to other PDEs, such as
convection-diffusion, solid mechanics or electromagnetism (Maxwell's
equations).

\CMS{Future work?}
}

\appendix

\section{2D Quadrilateral Under Uniaxial Stretch}\label{sec:stencil2d}
\begin{figure}[htb!]\label{fig:2dsa}
\centering
\includegraphics[trim=2 8 2 22,clip,scale=0.5]{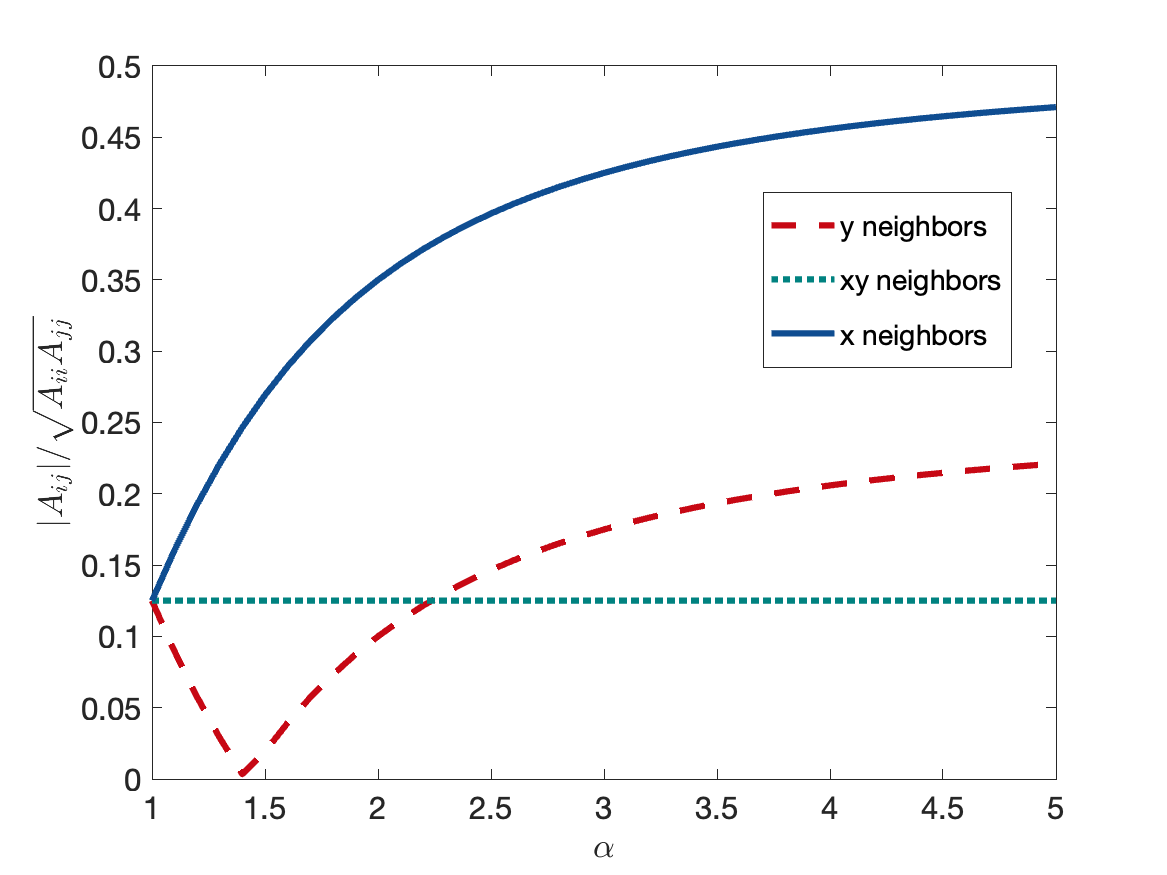}
\caption{Traditional SA criterion values for 2D stretched stencil. } 
\vskip -.1in
\end{figure}

Consider a 2D structured mesh where the $y$-direction is stretched by a constant factor $\alpha$, relative to the mesh spacing in the
$x$-direction, $h$. 
Discretization of a Poisson operator by first order nodal quadrilateral elements yields the element stiffness matrix
\begin{equation}
\frac{1}{6\alpha}
\left[ \begin {array}{cccc} 2\,{\alpha}^{2}+2&-2\,{\alpha}^{2}+1&-{
\alpha}^{2}-1&{\alpha}^{2}-2\\ \noalign{\medskip}-2\,{\alpha}^{2}+1&2
\,{\alpha}^{2}+2&{\alpha}^{2}-2&-{\alpha}^{2}-1\\ \noalign{\medskip}-{
\alpha}^{2}-1&{\alpha}^{2}-2&2\,{\alpha}^{2}+2&-2\,{\alpha}^{2}+1
\\ \noalign{\medskip}{\alpha}^{2}-2&-{\alpha}^{2}-1&-2\,{\alpha}^{2}+1
&2\,{\alpha}^{2}+2\end {array} \right]
\end{equation}
whose $(i,j)$ entries are defined by
$\int \grad\phi_j \cdot \grad\phi_i \;d\Omega_e$  where $\phi_k$ is the
$k^{th}$ vertex basis function, $\Omega_e$ is the stretched
element, and vertices are ordered in a
counter-clockwise fashion starting at the $(-1,-1)$ node on the
reference element.
This can be converted to an interior vertex stencil
by noting that
four elements contribute to the stencil central point while two elements contribute to the stencil sides,
leading to the stencil
\\[-12pt]
\begin{equation}\label{eq:stencil2d_A}
\frac{1}{6\alpha}
\StencilTwoD{8+8\alpha^2}{-4+2\alpha^2}{2-4\alpha^2}{-1-\alpha^2} .
\end{equation}
\vskip -1.001in
\hskip 1.49in
\hskip -.25in
{\begin{tikzpicture}
\draw[dashed,tealgreen,line width=2pt] (3.51,0.18) rectangle (5.10,0.83);
\draw[venetianred,dotted,line width=2pt] (6.36,0.51) node [ellipse,draw,minimum width=.7in,outer sep=-10pt, minimum height=.21in]{};
\draw[dashed,tealgreen,line width=2pt] (7.54,0.18) rectangle (9.10,0.83);
\draw[yaleblue,line width=2pt] (4.31,1.36) node [ellipse,draw,minimum width=.7in,outer sep=-10pt, minimum height=.21in]{};
\draw[yaleblue,line width=2pt] (8.37,1.36) node [ellipse,draw,minimum width=.7in,outer sep=-10pt, minimum height=.21in]{};
\draw[dashed,tealgreen,line width=2pt] (3.51,1.86) rectangle (5.10,2.51);
\draw[venetianred,dotted,line width=2pt] (6.31,2.20) node [ellipse,draw,minimum width=.7in,outer sep=-10pt, minimum height=.21in]{};
\draw[dashed,tealgreen,line width=2pt] (7.54,1.86) rectangle (9.10,2.51);
\end{tikzpicture}\label{fig:bleck2d}}
\\Here, identical stencil values are highlighted with identical dotted rectangles, dotted ovals, or solid ovals.
Figure~\ref{fig:2dsa}
plots the corresponding traditional SA criterion values. To label nearby $x$ neighbors as strong
and label distant $xy$ and $y$ neighbors as weak requires a threshold greater than $.25$ but less than $.5$
when $\alpha$ is large.  Unfortunately, this threshold would then label all off-diagonals as weak when $\alpha \approx 1$.

\begin{figure}\label{fig:live2dDistribLump}
\hskip -.35in \includegraphics[trim = 2.2in 0.9in 1.9in 0.6in, clip,scale=0.235]{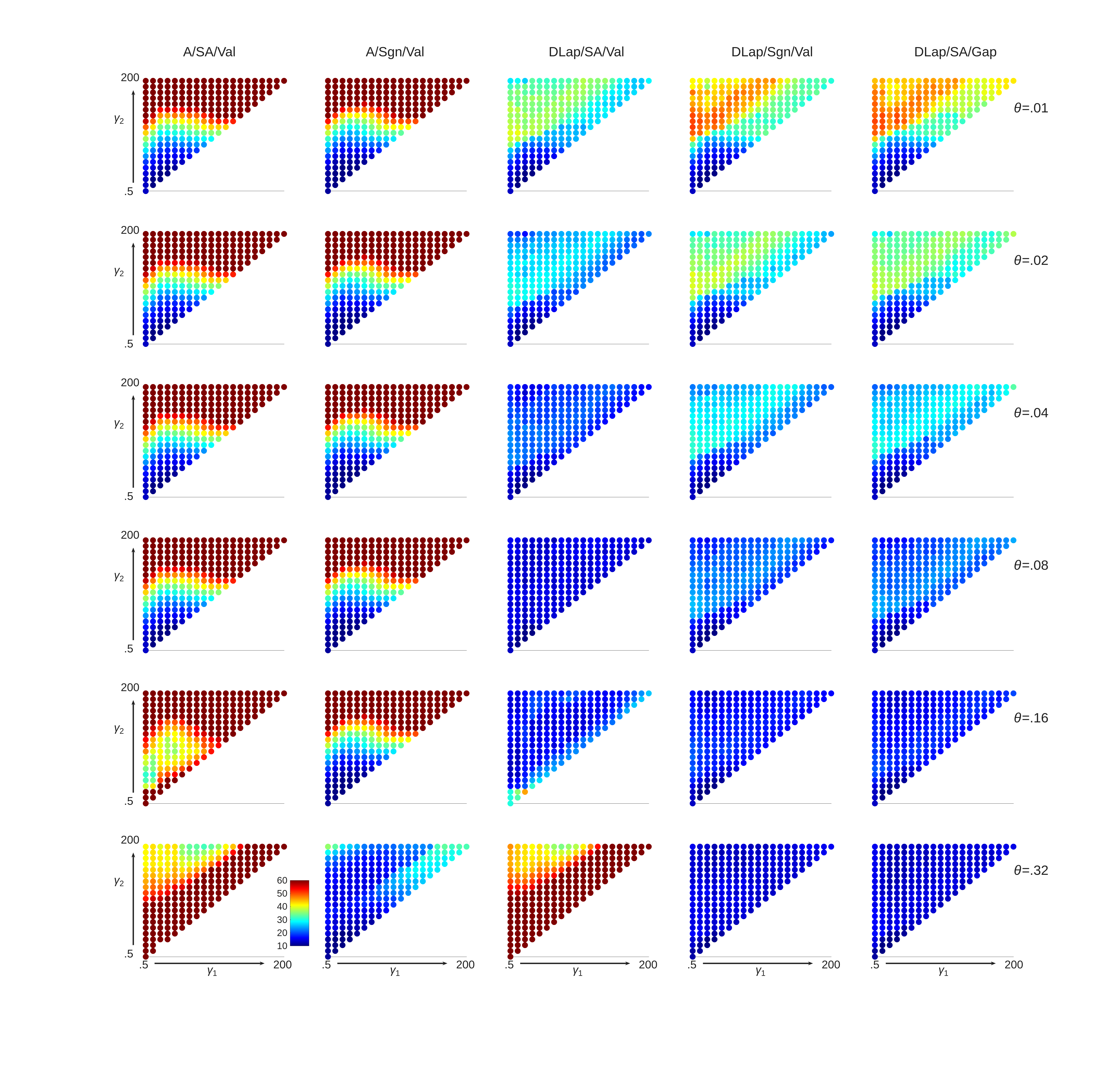}\vskip -.2in
\caption {Iterations{\hskip .01in}$\times${\hskip .01in}operator{\hskip-.02in} complexity, with the color scale from 10 to 60 for the 3D mesh
Columns represent algorithm choices (described in
Section~\ref{sec:results_2d}) and rows represent tolerances,
$\theta$. Distributed lumping is use for all experiments. See Section~\ref{sec:results_2d}.}
\end{figure}

Figure~\ref{fig:live2dDistribLump} depicts results using distributed lumping for the 2D stretched mesh and is the counterpart to
Figure~\ref{fig:live2d} which uses diagonal lumping. The two figures are similar, but one can identify cases
where distributed lumping improves the convergence (e.g., the speckled dots appearing for {\sf DLap/Sgn/Val} for $\theta = .16 ~\mbox{and}~ .32$
no longer appear speckled with distributed lumping).

\REMOVE{
Following \eqref{eq:entry_dlap}, the entries of distance
Laplacian matrix for an interior mesh vertex is
\REMOVE{
The off-diagonal
entries can be calculated directly, but the diagonal ``self'' edge
takes into account connectivity of other elements.
We note that each node connects to four
entries via mesh edges (two in the stretched direction and two in the
other direction) and four entries across the diagonal of the
element.
}
yields the stencil
\begin{equation}\label{eq:stencil2d_L}
h^{-2}
\StencilTwoD{\frac{2(\alpha^4+4\alpha^2+1)}{\alpha^2(1+\alpha^2)}}
{\alpha^{-2}}{1}{(1+\alpha^2)^{-1}}.
\end{equation}
We note that the ratio of off-diagonal elements to the diagonal is
independent of $h$, but does depend on $\alpha$.
}

\section{Element Stencil: 3D Hexahedron Under Uniaxial Stretch} \label{sec:stencil3d}

Consider a 3D structured mesh where the $z$-direction is stretched by a
constant factor $\alpha$ relative to the mesh spacing $h$ in the $x$- and $y$- directions.
Discretization of a Poisson operator by first order nodal finite elements on a hexahedral 
mesh yields
the element stiffness matrix 
\begin{equation}
\frac{h}{18\alpha}
\left[ \begin {array}{cccccccc} 4\,{\alpha}^{2}+2&-{\alpha}^{2}+1&-2
\,{\alpha}^{2}+{\frac{1}{2}}&-{\alpha}^{2}+1&2\,{\alpha}^{2}-2&-{
\frac {{\alpha}^{2}}{2}}-1&-{\alpha}^{2}-{\frac{1}{2}}&-{\frac {{
\alpha}^{2}}{2}}-1\\ \noalign{\medskip}-{\alpha}^{2}+1&4\,{\alpha}^{2}
+2&-{\alpha}^{2}+1&-2\,{\alpha}^{2}+{\frac{1}{2}}&-{\frac {{\alpha}^{2
}}{2}}-1&2\,{\alpha}^{2}-2&-{\frac {{\alpha}^{2}}{2}}-1&-{\alpha}^{2}-
{\frac{1}{2}}\\ \noalign{\medskip}-2\,{\alpha}^{2}+{\frac{1}{2}}&-{
\alpha}^{2}+1&4\,{\alpha}^{2}+2&-{\alpha}^{2}+1&-{\alpha}^{2}-{\frac{1
}{2}}&-{\frac {{\alpha}^{2}}{2}}-1&2\,{\alpha}^{2}-2&-{\frac {{\alpha}
^{2}}{2}}-1\\ \noalign{\medskip}-{\alpha}^{2}+1&-2\,{\alpha}^{2}+{
\frac{1}{2}}&-{\alpha}^{2}+1&4\,{\alpha}^{2}+2&-{\frac {{\alpha}^{2}}{
2}}-1&-{\alpha}^{2}-{\frac{1}{2}}&-{\frac {{\alpha}^{2}}{2}}-1&2\,{
\alpha}^{2}-2\\ \noalign{\medskip}2\,{\alpha}^{2}-2&-{\frac {{\alpha}^
{2}}{2}}-1&-{\alpha}^{2}-{\frac{1}{2}}&-{\frac {{\alpha}^{2}}{2}}-1&4
\,{\alpha}^{2}+2&-{\alpha}^{2}+1&-2\,{\alpha}^{2}+{\frac{1}{2}}&-{
\alpha}^{2}+1\\ \noalign{\medskip}-{\frac {{\alpha}^{2}}{2}}-1&2\,{
\alpha}^{2}-2&-{\frac {{\alpha}^{2}}{2}}-1&-{\alpha}^{2}-{\frac{1}{2}}
&-{\alpha}^{2}+1&4\,{\alpha}^{2}+2&-{\alpha}^{2}+1&-2\,{\alpha}^{2}+{
\frac{1}{2}}\\ \noalign{\medskip}-{\alpha}^{2}-{\frac{1}{2}}&-{\frac {
{\alpha}^{2}}{2}}-1&2\,{\alpha}^{2}-2&-{\frac {{\alpha}^{2}}{2}}-1&-2
\,{\alpha}^{2}+{\frac{1}{2}}&-{\alpha}^{2}+1&4\,{\alpha}^{2}+2&-{
\alpha}^{2}+1\\ \noalign{\medskip}-{\frac {{\alpha}^{2}}{2}}-1&-{
\alpha}^{2}-{\frac{1}{2}}&-{\frac {{\alpha}^{2}}{2}}-1&2\,{\alpha}^{2}
-2&-{\alpha}^{2}+1&-2\,{\alpha}^{2}+{\frac{1}{2}}&-{\alpha}^{2}+1&4\,{
\alpha}^{2}+2\end {array} \right]
\end{equation}
where $\Omega_e$ is the stretched element and 
degrees of freedom are ordered in a front-to-back, counter-clockwise fashion starting at the $(-1,-1,-1)$ node on the
reference element.  
This can be converted into an interior stencil by viewing the stencil as a cube and 
noting that stencil central point has eight element contributions, the six stencil values associated with 
face centers have two element contributions, the twelve stencil values associated with edge centers have four element
contributions, while the eight stencil values associated with corners have one element contribution.
This yields the stencil depicted in Figure~\ref{fig:bleez}. 
\section{Pamgen Mesh Generation Templates}

\subsection{2D Stretched Brick Mesh}\label{sec:pamgen2d}
Section~\ref{sec:results_2d} meshes are generated with the following
Pamgen \cite{Pamgen} template:
\begin{verbatim}
mesh
  brick
  numx 3
    xblock 1 1.0, interval 10
    xblock 2 {3.0*(GAMMA1+1)}, first size .1, last size {GAMMA1/10}
    xblock 3 {GAMMA1}, interval 10
  numy 3
    yblock 1 1.0, interval 10
    yblock 2 {3.0*(GAMMA2+1)}, first size .1, last size {GAMMA2/10}
    yblock 3 {GAMMA2}, interval 10
  end
  set assign
    sideset, jlo, 2
  end
end
\end{verbatim}
Braces denote preprocessor substitution and the variables \texttt{GAMMA1} and
\texttt{GAMMA2} are varied logarithmically between 0.5 and 200.

\subsection{3D Stretched Brick Mesh}\label{sec:pamgen3d}
Section~\ref{sec:results_3d} meshes are generated with the following Pamgen template:
\begin{verbatim}
mesh
  brick
  numx 3
    xblock 1 1.0, interval 10
    xblock 2 {3.0*(GAMMA1+1)}, first size .1, last size {GAMMA1/10}
    xblock 3 {GAMMA1}, interval 10
  numy 3
    yblock 1 1.0, interval 10
    yblock 2 {3.0*(GAMMA2+1)}, first size .1, last size {GAMMA2/10}
    yblock 3 {GAMMA2}, interval 10
  numz 1
    zblock 1 8.0, interval 80
  end
  set assign
    sideset, jlo, 2
  end
end
\end{verbatim}

\subsection{3D Radial Trisection Mesh}\label{sec:pamgen3d_radtri}
Section~\ref{sec:results_3d_radtri} meshes are generated with the following Pamgen template:
\begin{verbatim}
  mesh
    radial trisection
      trisection blocks, {GAMMA2}
      transition radius, 6.
      numz 1
        zblock 1 20. interval 20
      numr 2
        rblock 1 8.0 interval 4
        rblock 2 {200*floor((1+GAMMA1)/2)} first size 1.0 last size {GAMMA1}
      numa 1
        ablock 1 90. interval 24
    end
    set assign
      sideset, ihi, 4
    end
  end
\end{verbatim}
The variable \texttt{GAMMA2} indicates the number of trisection blocks.  The variable \texttt{GAMMA1}
is varied logarithmically up to 100, depending on the number of
trisection blocks (this is due to limits of Pamgen).  For one, two, three,
four and six blocks, \texttt{GAMMA1} correspondingly begins at 1.0658, 1.1875, 1.3231,
1.4742 and 1.6426.


\bibliographystyle{siamplain}
\bibliography{references,newrefs}
\end{document}